\documentclass[12pt]{article}
\usepackage{amsmath, latexsym, amssymb, amsthm}
\usepackage[usenames]{color}
\usepackage[normalem]{ulem}
\usepackage[all]{xy}

\newcommand{\vp}{\varphi}

\newcommand\ot{\otimes}

\newcommand{\wh}{\widehat}

\newcommand{\Aut}{\operatorname{Aut}}
\newcommand{\Res}{\operatorname{Res}}
\newcommand{\infl}{\operatorname{Infl}}
\newcommand{\id}{\operatorname{id}}
\newcommand{\Id}{\operatorname{Id}}

\newcommand{\bean}{\begin{eqnarray}}
\newcommand{\eean}{\end{eqnarray}}
\newcommand{\bea}{\begin{eqnarray*}}
\newcommand{\eea}{\end{eqnarray*}}
\newcommand{\bsa}{\begin{subarray}{c}}
\newcommand{\esa}{\end{subarray}}
\newcommand{\bi}{\begin{itemize}}
\newcommand{\ei}{\end{itemize}}
\newcommand{\sr}{\stackrel}

\newtheorem{thm}{Theorem}[section]
\newtheorem{lem}[thm]{Lemma}
\newtheorem{prop}[thm]{Proposition}

\newtheorem{rem}[thm]{Remark}
\newtheorem{cor}[thm]{Corollary}

\renewcommand{\a}{\alpha}
\newcommand{\w}{\omega}
\newcommand{\g}{\gamma}
\newcommand{\e}{\epsilon}

\newcommand\ol[1]{\overline{#1}}
\newcommand\C[1]{#1\mbox{\rm-\bf{Mod}}}
\newcommand\Cf[1]{{\C{#1}_\mathsf{fin}}}
\newcommand\Co[1]{{}^#1\!\mathcal{M}}
\newcommand\QH[2]{{\BC^{#1}\!{#2}}}
\newcommand{\ve}{\varepsilon}
\newcommand\la{\widetilde{\leftharpoonup}}

\newcommand{\inv}{^{-1}}
\newcommand{\du}{^\vee}
\newcommand{\bidu}{^{\vee\vee}}

\renewcommand{\to}{\rightarrow}
\newcommand\rt{\rtimes}

\newcommand{\BC}{\mathbb{C}}
\newcommand{\BF}{\mathbb{F}}
\newcommand{\BZ}{\mathbb{Z}}
\newcommand\CC{\mathcal{C}}
\newcommand\DD{\mathcal{D}}
\newcommand\FF{\mathcal{F}}

\DeclareMathOperator\MM{\mathcal{M}}
\DeclareMathOperator\ev{{\operatorname{ev}}}
\DeclareMathOperator\db{\operatorname{db}}

\def\ssp{\def\baselinestretch{1.0}\large\normalsize}

\title{\sc On the Gauge Equivalence of  Twisted Quantum Doubles of
Elementary Abelian and  Extra-Special $2$-Groups}

\author{\sc Christopher Goff\thanks{Supported by an Eberhardt
Research Fellowship from the University of the Pacific.} \\
Department of Mathematics, University of the Pacific \\
    \sc Geoffrey Mason\thanks{Supported by a NSF grant and faculty
research funds granted by the University of
California at Santa Cruz.}  \\
    Department of Mathematics, UC Santa Cruz \\
    \sc  Siu-Hung Ng\thanks{Supported by NSA H98230-05-1-0020}  \\
      Department of Mathematics, Iowa State University}
\date{}
\begin{document}
\ssp \maketitle

\begin{abstract}
      \noindent
      We establish braided tensor equivalences among module categories over
      the  twisted quantum double of  a finite group defined by an
   extension of a group  $\ol{G}$ by an abelian group,
      with  $3$-cocycle inflated from a $3$-cocycle on $\ol{G}$. We
also prove that the canonical
      ribbon structure of the module category of any twisted quantum
double of a finite group is preserved by
      braided tensor equivalences. We give two main applications:
first,  if $G$ is an extra-special $2$-group
      of width at least $2$, we show that the quantum double of $G$ twisted by
      a $3$-cocycle $\omega$ is gauge equivalent to a twisted quantum
double of an elementary abelian $2$-group if, and only if,
$\omega^2$ is trivial; second,  we discuss the gauge equivalence
classes of twisted quantum doubles of groups of order $8$, and
classify the braided tensor equivalence classes of these
quasi-triangular quasi-bialgebras. It turns out that there are
exactly $20$ such
   equivalence classes.
    \end{abstract}

\section{ Introduction}

Given a finite group $G$ and a (normalized) $3$-cocycle $\omega \in
Z^3(G, \mathbb{C}^*)$, there is associated the \emph{twisted quantum
double} $D^{\omega}(G)$. This is a certain braided quasi-Hopf
algebra introduced by Dijkgraaf-Pasquier-Roche in the context of
orbifold conformal field theory [DPR]. For a holomorphic vertex
operator algebra (VOA)  $V$ admitting a faithful action of $G$ as
automorphisms, one expects that the orbifold VOA $V^G$ has a module
category $\C{V^G}$ which (among other things) carries the structure
of a braided tensor category, and that there is an equivalence of
braided tensor categories $\C{V^G} \cong \C{D^{\omega}(G)}$ for some
choice of $\omega$. The conjectured equivalence of categories is
deep, and motivates the results of the present paper. These,
however, are concerned purely with twisted quantum doubles and their
module categories. No familiarity with VOA theory is required to
understand our main results, and we use the language of VOAs purely
for motivation. For further background on the connection between
VOAs and quasi-Hopf algebras, see [DPR],
[M1], [MN1] and [DM].\\

There is an interesting phenomenon, akin to mirror symmetry, that
arises as follows:
    we are given two pairs $(V, G), (W, H)$, where $V, W$ are
holomorphic VOAs with finite,
    faithful automorphism groups $G, H$ respectively, together with an
isomorphism
    of VOAs,  $V^G \cong W^H$.  If the conjectured equivalence of
orbifold categories is true, there
    must also be a tensor equivalence
    \begin{eqnarray}\label{tensequiv}
    \C{D^{\omega}(G)} \cong \C{D^{\eta}(H)}
    \end{eqnarray}
for some choice of $3$-cocycles $\omega, \eta$ on $G$ and $H$
respectively. Conversely, deciding when equivalences such as
(\ref{tensequiv}) can hold gives information about the VOAs and the
cocycles that  they determine. This is an interesting problem in its
own right, and is the one we consider here.\\

    The case in which the two twisted doubles in question are
\emph{commutative} was treated at length in [MN1]. Here we are
concerned with a particular case of the more difficult situation in
which the twisted doubles are not commutative. It arises in orbifold
theory when one applies the $\BZ_2$-orbifold construction to a
holomorphic lattice theory $V_L$ ($L$ is a self-dual even lattice)
[FLM], [DGM].  One takes $G$ to be the \emph{elementary abelian
$2$-group} $L/2L \oplus \BZ_2$, and it turns out that $H$ is an
\emph{extra-special} $2$-group. We study this situation in the
present paper. Changing notation,  one of our main results (Theorems
\ref{exspecqdbleisom} and \ref{FSexpcomp}) is as follows:
\begin{equation}\label{introthm5.3}
\parbox[c]{5.2in}{Let $Q$ be an extra-special group of order $2^{2l+1}$
\emph{not} isomorphic to the dihedral group of order $8$, and let
$\eta$ be \emph{any} 3-cocycle such that ${[\eta]} \in H^3(Q,
\mathbb{C}^*)$ has order at most 2. Then there are $\mu \in Z^3(E,
\mathbb{C}^*)$, $E = \BZ_2^{2l+1}$, and a braided tensor equivalence
$\C{D^{\mu}(E)}  \cong \C{D^{\eta}(Q)}$. Such a tensor equivalence
does  \emph{not} exist if $[\eta]$ has order \emph{greater} than
$2$.}
\end{equation}

As is well-known (cf. \cite{EG}, \cite{NS2}), such an equivalence of
braided tensor categories corresponds to a \emph{gauge equivalence}
of the twisted doubles as quasi-triangular quasi-bialgebras. The
proof of (\ref{introthm5.3}) relies on Theorem \ref{thmfamily}
together with the gauge-invariance of \emph{Frobenius-Schur
exponents} (cf. \cite{NS3}) of semisimple quasi-Hopf algebras. In
Theorem \ref{thmfamily} we establish the existence of gauge
equivalences of quasi-triangular quasi-bialgebras among twisted
quantum doubles
   $D^{\omega}(G)$
where $G$ is a group defined as an extension of ${\ol G}$ by an
abelian group and the $3$-cocycle $\omega \in Z^3(G, \BC^*)$ is
inflated from a 3-cocycle on $\ol G$. Several results of Schauenburg
(\cite{Schauenburg01}, \cite{SchAdv02})  play a significant r\^{o}le
in the proof.\\

As long as $l \geq 2$, the group of elements of $H^3(Q,
\mathbb{C}^*)$ of order at most $2$ has index $2$ in the full degree
$3$ cohomology ([HK]). The proof of the last assertion of
(\ref{introthm5.3}) involves some computations involving the
cohomology of  $Q$. We use these to calculate (Theorem
\ref{FSexpcomp}) the Frobenius-Schur exponents of $D^{\omega}(G)$
where $G$ is either extra-special group or elementary abelian, and
$\omega$ is any  $3$-cocycle.
    Frobenius-Schur
indicators, their higher analogs, and  Frobenius-Schur exponents
have been investigated in \cite{MN2}, \cite{NS1}, \cite{NS2},
\cite{NS3} in the general context of semisimple quasi-Hopf algebras
and pivotal fusion categories. They provide valuable gauge
invariants which are
    reasonably accessible in the case of twisted quantum doubles.\\

The case $l=1$ is exceptional in several ways,  and we consider the
    problem of enumerating the gauge equivalence classes defined by
quasi-Hopf algebras obtained
    by twisting the quantum doubles
of $\BZ_2^3, Q_8$ and $D_8$.  There are $88$ such quasi-Hopf
algebras which are noncommutative, corresponding to the $64$
nonabelian cohomology classes for $\BZ_2^3$,
    together with $8$ classes for $Q_8$ and $16$ for $D_8$.
     Some of the subtlety of this problem, which remains open,
can be illustrated by observing that \emph{all} $88$  twisted
doubles have the \emph{same} fusion algebra, moreover some of them
are isomorphic as \emph{bialgebras} yet are not gauge equivalent. We
will see that there are at least $8$, and no more than $20$, gauge
equivalence classes. This makes use of (\ref{introthm5.3}) to
identify equivalence classes of quasi-bialgebras, together with
Frobenius-Schur indicators and their higher analogs (loc. cit) to
distinguish between equivalence classes. Now there is a
\emph{canonical braiding} of these quasi-Hopf algebras, and using
the \emph{invariance} of the canonical ribbon structure under
braided  tensor equivalences, which we establish separately,  we
show
   that the 20 gauge equivalence classes constitute a complete
list of gauge equivalence classes of the quasi-triangular
quasi-bialgebras under consideration.

\medskip
The paper is organized as follows: following some background, in
Section \ref{sectionGE} we give the proof of  Theorem
\ref{thmfamily}.
   There are connections here to some results of Natale \cite{N}.
In Section \ref{sectionBI} we discuss a variant of Theorem
\ref{thmfamily} involving bialgebra isomorphisms. In Section
\ref{sectionES} we give the proofs of  Theorems
\ref{exspecqdbleisom} and  \ref{FSexpcomp}. The question of the
gauge equivalence classes for twisted doubles of groups of order $8$
is presented in Section \ref{section64}.
    In Section \ref{ribbon}
we show that the twist associated with a pivotal braided tensor
category is preserved by tensor equivalences which preserve both
pivotal structure and braiding. In particular, the \emph{canonical
ribbon structure} of the module category of a semisimple
quasi-triangular quasi-Hopf algebra is preserved by braided tensor
equivalences. In an Appendix we give a complete set of
Frobenius-Schur indicators, their higher analogs, and the scalars of
the ribbon structures for the particular quasi-Hopf algebras
    arising from twisted doubles of groups of order $8$.
This data is interesting in its own right, and may be useful in the
future.

    \medskip
    We take this opportunity to thank the referee for several useful
comments on a prior version
    (math.QA/0603191)  of
    the current paper. In particular, the referee pointed out the
possibility of utilizing results
    of Schauenburg to significantly improve on our earlier version of
Theorem 2.1.

    \medskip
     We use the following conventions: all linear spaces are defined over the
complex numbers $\mathbb{C}$ and have finite dimension; all
     groups are finite;
     all cocycles for a group $G$ are defined with respect to a
\emph{trivial} $G$-module
     and are \emph{normalized}; a \emph{tensor category} is a
$\BC$-linear monoidal category; a \emph{tensor functor} is a
$\BC$-linear monoidal functor; braided quasi-Hopf algebra means
quasi-triangular quasi-Hopf algebra. We generally do not
differentiate between a cocycle and the cohomology class it defines,
and often use the isomorphism
     $H^{n+1}(G, \BZ) \cong H^{n}(G, \mathbb{C}^*)$ without
explicit reference to it.
     Unexplained notation is as in [K].

\section{A Family of Gauge Equivalences}\label{sectionGE}

We first review some facts about twisted quantum doubles ([DPR],
[K]).
    Fix a group $G$ and $3$-cocycle $\omega \in
Z^3(G, \mathbb{C}^*)$.
    $D^{\omega}(G)$ is a crossed product with underlying linear space
    $\mathbb{C}G^* \otimes \mathbb{C}G$,
    where the dual group algebra
$\mathbb{C}G^*$ has basis $e(g)$ dual to the group elements $g \in
G$. Multiplication, comultiplication, associator, counit, antipode,
$R$-matrix, $\alpha$ and $\beta$ are given by the following
formulas:
\begin{eqnarray*}
&&e(g) \otimes x. e(h) \otimes y = \theta_g(x, y)
\delta_g^{xhx^{-1}}
e(g) \otimes xy, \\
&&\Delta (e(g) \otimes x) = \sum_{\bsa h, k  \\ hk = g \esa}
\gamma_x(h, k) e(h)
\otimes x \otimes e(k) \otimes x, \\
&&\Phi = \sum_{g, h, k} \omega(g, h, k)^{-1} e(g) \otimes 1 \otimes
e(h) \otimes 1 \otimes e(k) \otimes 1, \\
&&\epsilon(e(g) \otimes x) = \delta_g^{1}, \\
&&S(e(g) \otimes x) = \theta_{g^{-1}}(x, x^{-1})^{-1} \gamma_x(g,
g^{-1})^{-1} e(x^{-1}g^{-1}x)
\otimes x^{-1},\\
&&R = \sum_{g, h} e(g) \otimes 1 \otimes e(h) \otimes g, \\
&& \alpha = 1, \ \beta = \sum_g \omega(g, g^{-1}, g) e(g) \otimes 1.
\end{eqnarray*}
Here, $\theta_g, \gamma_g$ are $2$-cochains on $G$ derived from
$\omega$ via
\begin{eqnarray}
&& \theta_g(x,y) = \frac{\omega(g,x,y)\omega(x,y,(xy)^{-1}gxy)}
{\omega(x,x^{-1}gx,y)},
\label{thetaformula} \\
&&\gamma_g(x, y) = \frac{\omega(x, y, g) \omega(g, g^{-1}xg,
g^{-1}yg)}{\omega(x, g, g^{-1}yg)}. \label{gammaformula}
\end{eqnarray}

For any quasi-Hopf algebra $H = (H, \Delta, \epsilon, \Phi, \alpha,
\beta, S)$ and any unit $u \in H$, we can \emph{twist} $H$ by $u$ to
obtain a second quasi-Hopf algebra $H_u =  (H, \Delta, \epsilon,
\Phi, u\alpha, \beta u^{-1}, S_u),$ where $S_u(h) = uS(h)u^{-1}$ for
$h \in H.$ See [D] for details. We obtain another quasi-Hopf algebra
$H_F = (H, \Delta_F, \epsilon, \Phi_F, \alpha_F, \beta_F, S)$ via a
\emph{gauge transformation} determined by suitable unit $F \in H
\otimes H$ ([K]).\\

The notion dual to quasi-bialgebra is coquasi-bialgebra (cf.
\cite{Maji93}). A coquasi-bialgebra  $H=(H, \phi)$ is a bialgebra
$H$ (although the multiplication is not necessarily associative),
together with a convolution-invertible map $\phi: H^{\ot 3} \to
\BC$, called \emph{coassociator}, which satisfies the dual
conditions of associators. The category $\Co{H}$ of left comodules
of a coquasi-bialgebra $H$ is a tensor category with associativity
isomorphism determined by the coassociator $\phi$. Suppose that $K$
is a Hopf algebra, $\iota: K \to H$ a bialgebra map such that $\phi
\circ (\iota \ot \id_H \ot \id_H )=\ve_K \ot \ve \ot \ve$ ($\ve$ is
the counit of $H$), and that there exists a convolution invertible
left $K$-module map $\pi: H \to K$ satisfying $\pi(1)=1$, $\ve \pi =
\ve$ and
$$
\phi(g \ot x \ot h_{(1)})\pi(h_{(2)}) =\phi(g \ot x \ot
\pi(h)_{(1)})\pi(h)_{(2)}, \quad\text{for } g,h\in H, \, x \in K.
$$
Then by \cite[Corollary 3.4.4]{SchAdv02}, $\sideset{_K^H}{_K}\MM$ is
a tensor category  and there exists a coquasi-bialgebra $B$ such
that $\Co{B}$ is equivalent  as tensor category to
$\sideset{_K^H}{_K}\MM$. The structure maps of $B$ are given in
\cite[3.4.2 and 3.4.5]{SchAdv02}. This result of Schauenburg plays
an
    important r\^{o}le in the proof of
Theorem \ref{thmfamily}.\\

For a group  $G$, and $3$-cocycle $\w \in Z^3(G, \BC^*)$, the
bialgebra $\BC G$ together with the \emph{coassociator}
\begin{equation}\label{eq:coassociator}
\phi = \sum_{x,y, z} \w\inv(x,y,z) e(x) \ot e(y)\ot e(z) \in (\BC
G^{\ot 3})^*
\end{equation}
is a \emph{coquasi-bialgebra}. We simply write $\QH{\w}{G}$ for this
coquasi-bialgebra. Note that $\Co{\QH{\w}{G}}$ is tensor equivalent
to the category of modules over
    the quasi-bialgebra
$\BC G^*$ equipped with the associator $\phi$, and every
$\QH{\w}{G}$-comodule is a $G$-graded vector space. The
associativity isomorphism $\Phi: (U \ot V) \ot W \to U \ot (V \ot
W)$ is given by
$$
\Phi(u \ot v \ot w) = \w\inv(g, g', g'') u \ot v \ot w
$$
for homogeneous elements $u \in U$, $v \in V$, $w \in W$ of degree
$g, g', g''$ respectively. Moreover, the center of $\Co{\QH{\w}{G}}$
is equivalent to
$\C{D^\w(G)}$ as braided tensor categories.\\

Let ${\ol{G}}$ be a group, $A$ a right ${\ol{G}}$-module with
${\ol{G}}$-action $\triangleleft$, and $E$ the semidirect product $A
\rtimes {\ol{G}}$ with underlying set  $ \ol G \times A$ and
multiplication
$$
(x,a)(y, b) = (xy, (a \triangleleft y)b)
$$
for $x,y \in {\ol{G}}, a,b \in A$. Then $E$ fits into a \emph{split}
exact sequence of groups:
\begin{equation}\label{eq:e1}
1 \longrightarrow A \longrightarrow E \longrightarrow {\ol{G}}
\longrightarrow 1\,.
\end{equation}

The character group $\hat{A} = \mbox{Hom}(A, \mathbb{C}^*)$ of $A$
admits a left ${\ol{G}}$-module structure $\triangleright$ given by
$$
(h \triangleright \chi)(a) = \chi(a \triangleleft h)\,.
$$
Let $G$ be an extension of ${\ol{G}}$ by $\hat{A}$ associated with
the exact sequence of groups:
\begin{equation}\label{eq:e2}
1 \longrightarrow \hat{A} \longrightarrow G \longrightarrow {\ol{G}}
\longrightarrow 1\,.
\end{equation}
Display (\ref{eq:e2})  determines  a $2$-cocycle $\e \in
Z^2({\ol{G}}, \hat{A})$, i.e. a function $\epsilon: {\ol{G}} \times
{\ol{G}} \to \hat{A}$ satisfying:
\begin{gather}
\epsilon(x, 1_{\ol{G}})=\epsilon(1_{\ol{G}}, y)=1_{\hat{A}}\,,\\
\label{cocycle} (x \triangleright \epsilon(y,z)) \epsilon(x,yz)
=\epsilon(xy,z)\epsilon(x,y)
\end{gather}
for all   $x,y,z \in {\ol{G}}$. Then  $G$ is isomorphic to a group
$\hat{A} \rtimes_{\epsilon} {\ol{G}}$ which has $\hat{A}\times \ol
G$ as underlying set
    and multiplication
$$
(\mu , x)(\nu, y) = (\mu (x\triangleright \nu)\epsilon(x,y), xy)
$$
for $x, y \in {\ol{G}}, \mu, \nu \in \hat{A}$. Moreover, the exact
sequence of groups
\begin{equation}\label{eq:e3}
1 \longrightarrow \hat{A} \sr{\iota}{\longrightarrow}
\hat{A}\rtimes_\epsilon {\ol{G}} \sr{ p}{\longrightarrow} {\ol{G}}
\longrightarrow 1
\end{equation}
is equivalent to \eqref{eq:e2}, where $\iota: \hat{A} \to
\hat{A}\rtimes_\epsilon {\ol{G}}$ and $p: \hat{A}\rtimes_\epsilon
{\ol{G}} \to {\ol{G}}$ are given by
$$
\iota(\mu) = (\mu, 1), \quad\text{and}\quad p(\mu, x)=x
$$
for $\mu \in \hat{A}, x \in {\ol{G}}$.\\

     We are going to compare twisted quantum doubles
of $G=\hat{A} \rtimes_{\epsilon} {\ol{G}}$ and $E=A \rtimes
{\ol{G}}$. Suppose that $\zeta \in Z^3({\ol{G}}, \mathbb{C}^*)$. We
may \emph{inflate} $\zeta$ along the projections $E \rightarrow
{\ol{G}}, G \rightarrow {\ol{G}}$ in (\ref{eq:e1}) and (\ref{eq:e3})
to obtain $3$-cocycles
    \begin{eqnarray*}
\zeta_E &=& \infl_{\ol{G}}^E \zeta \in Z^3(E, \mathbb{C}^*), \\
     \zeta_G &=& \infl_{\ol{G}}^G \zeta \in Z^3(G, \mathbb{C}^*).
\end{eqnarray*}
One can also check directly that the $3$-cochain $\omega$ on $E$
defined by
\begin{eqnarray}\label{cocycledef}
\w((h_1,a_1),(h_2,a_2),(h_3,a_3)) = \epsilon(h_2,h_3)(a_1),  \ h_i
\in {\ol{G}}, a_i \in A,
\end{eqnarray}
    actually lies in $Z^3(E, \mathbb{C}^*)$. Moreover, for any $\sigma
\in \Aut_{\ol{G}}(A)$,
    $$\sigma \w((h_1,a_1),(h_2,a_2),(h_3,a_3)) =
\epsilon(h_2,h_3)(\sigma(a_1))$$ also defines  a normalized
    3-cocycle on $E$, and we let $\Omega=\{\sigma \w \mid \sigma \in
\Aut_{\ol{G}}(A)\}$.\\

We can now state our first main result. Notation and assumptions are
as above.
\begin{thm}\label{thmfamily}
    Let $\zeta \in Z^3({\ol{G}}, \mathbb{C}^*)$  and  $\ol\w \in \Omega$.
    Then the braided tensor categories $\C{D^{\zeta_G}(G)}$  and
$\C{D^{\ol\w\zeta_E}(E)}$ are equivalent.
\end{thm}
\begin{proof} Bearing in mind the definition $G=\hat{A}
\rtimes_\epsilon {\ol{G}}$, it follows
     from the definition of inflated cocycle that
     $$
\zeta_G((\mu, x), (\mu', x'), (\mu'', x''))=\zeta(x,x',x'')
$$
for $(\mu, x), (\mu', x'), (\mu'', x'') \in G$.

     \medskip
Consider the coquasi-bialgebra $H=\QH{\zeta_G}{G}$ with coassociator
$\phi$ as in \eqref{eq:coassociator}, and the Hopf algebra $K=\BC
\hat{A}$. Let $\sigma \in \Aut_{\ol{G}}(A)$. We define $\BC$-linear
maps $\iota_\sigma: K \to H$ by $\iota_\sigma(\mu) = (\mu \circ
\sigma\inv, 1)$ for $\mu \in \hat{A}$, and $\pi_\sigma: H \to K$ by
$\pi_\sigma(\mu, x) = \mu\circ \sigma$ for $(\mu, x) \in G$.  It is
easy to see that $\iota_\sigma$ is a bialgebra map which satisfies
$$
\phi\circ (\iota_\sigma \ot \id_H \ot \id_H)= \ve_K \ot \ve \ot
\ve\,.
$$
Furthermore, $\pi_\sigma$ is a $K$-module coalgebra map and is
convolution invertible with inverse $\bar \pi_\sigma: H \to K$ given
by $\bar \pi_\sigma(\mu, x) =\mu\inv\circ \sigma$. Clearly,
$\pi_\sigma(1_H)=1_K$ and $\ve_K \circ \pi_\sigma = \ve_H$. By the
definitions of $\phi$ and $\zeta_G$, for $g, h \in G$ and $\mu \in
\hat{A}$, we have
$$
\phi(g \ot (\mu, 1) \ot h)\pi_\sigma (h) =\pi_\sigma(h) = \phi(g \ot
(\mu, 1) \ot \pi_\sigma(h))\pi_\sigma(h)\,.
$$
It follows from \cite[Corollary 3.4.4]{SchAdv02} that there exists a
coquasi-bialgebra $B$ such that
$$
\sideset{_K^H}{_K}\MM \cong \Co{B}
$$
as tensor categories. By \cite{Schauenburg01}, we also have the
equivalence
$$
\mathcal{Z}(\sideset{_K^H}{_K}\MM) \cong \mathcal{Z}(\Co{H})
$$
of braided tensor categories, where $\mathcal{Z}(\CC)$ denotes the
center of the monoidal category $\CC$.
   Since $\C{D^{\zeta_G}(G)}$ and $\mathcal{Z}(\Co{H})$
are equivalent braided tensor categories (cf. \cite{Maji98}), we
have
$$
\C{D^{\zeta_G}(G)} \cong \mathcal{Z}(\Co{B})
$$
as braided tensor categories.

     \medskip
We proceed to show that $B \cong \QH{\w'}{E}$ using \cite[3.4.2 and
3.4.5]{SchAdv02},  where $\w'=\sigma\tau\w\cdot \zeta_E$ and $\tau:
A \to A$, $\tau: a \mapsto a\inv$. We will use the notation defined
in \cite{SchAdv02} for the remaining discussion. By \cite[Lemma
3.4.2]{SchAdv02}, $B=\BC {\ol{G}} \ot K^*$ as vector space, and we
write $x \rtimes a$ for $x \ot a$ whenever $x \in {\ol{G}}$ and $a
\in A$. Here, we have used the canonical identification of the Hopf
algebra $K^*$ with $\BC A$. Note that $\BC {\ol{G}} \cong H/K^+ H$
as coalgebras. Since $\iota_\sigma(\hat{A})$ is a normal subgroup of
$G$, $K^+H = HK^+$. Thus, the right $K$-action on $H/K^+ H$ is
trivial, and so the associated left $K^*$-comodule structure on $\BC
{\ol{G}}$ is also trivial. Hence, $\tilde \rho(q) = 1_A \ot q$ for
$q \in \BC {\ol{G}}$. For any $x \in {\ol{G}}$, $\mu, \nu \in
\hat{A}$,
$$
\phi((1, x)\ot (\mu, 1)\ot (\nu, 1)) = 1\,.
$$
So the $\tilde{\a} : \BC {\ol{G}} \to K^* \ot K^*$ defined in
\cite[Lemma 3.4.2]{SchAdv02} for our context is given by $\tilde \a
: x \mapsto 1 \ot 1$. Thus, the comultiplication $\tilde \Delta$ and
counit $\tilde \ve$ of $B$ are given by
$$
\tilde \Delta(x \rtimes a) = x \rtimes a \ot x \rtimes a \quad
\text{and}\quad \tilde\ve(x \rtimes a) = 1
$$
for $x \in {\ol{G}}$ and $a \in A$.  With the \emph{cleaving} map
$\pi_\sigma$, the associated map $j_\sigma: \BC {\ol{G}} \to H$ is
given by $j_\sigma(x) = (1,x)$ for $x \in {\ol{G}}$, and so
$$
x \rightharpoonup \mu = \pi_\sigma(j_\sigma(x)(\mu\circ \sigma, 1))=
(x \triangleright (\mu \circ \sigma))\circ \sigma\inv = x
\triangleright \mu
$$
for $\mu \in \hat{A}$. Thus for $a \in A$, $x \in {\ol{G}}$ and $\mu
\in \hat{A}$, we have
$$
\mu(a \la x)=(x \rightharpoonup \mu)(a) = (x \triangleright \mu)(a)
= \mu(a \triangleleft x)
$$
and hence
$$
a \la x = a \triangleleft x\,.
$$
Recall that
$$
e(a) = \frac{1}{|A|}\sum_{\chi \in \hat{A}}
\frac{1}{\chi(a)}\chi\,,\quad a \in A,
$$
form a dual basis of $A$. It follows from \cite[Theorem
3.4.5]{SchAdv02} that the multiplication on $B$ is given by
\begin{eqnarray*}
(x \rtimes a)\cdot  (y \rtimes b) & = & \sum_{ c\in A}
\ol{j_\sigma(x)j_\sigma(y)} \rtimes c (a \widetilde{\leftharpoonup}
y)b \phi(j_\sigma(x) \ot j_\sigma(y) \ot \iota_\sigma(e(c)) )
\pi_\sigma(j_\sigma(y))(a)\\
   & = &
\frac{1}{|A|}\sum_{ c\in A} xy \rtimes c (a \triangleleft y)b
\sum_{\chi \in \hat{A}} \frac{1}{\chi(c)}\phi((1,x) \ot (1,y) \ot
(\chi\circ \sigma\inv, 1) )\\
& = &
\frac{1}{|A|}\sum_{ c\in A} xy \rtimes c (a \triangleleft y)b |A|\delta_c^1\\
&=& xy \rtimes (a \triangleleft y)b,
\end{eqnarray*}
and the coassociator is
\begin{eqnarray*}
\tilde\phi (x \rt a \ot y \rt b \ot z \rt c) & = & \phi(j_\sigma(x)
\ot j_\sigma(y) \ot
j_\sigma(z))\pi_\sigma(j_\sigma(y)j_\sigma(z))(a)\pi_\sigma(j_\sigma(z))(b)
\ve(c)\\
& = & \zeta_G((1,x), (1,y), (1,z))\inv \pi_\sigma(\e(y, z), yz)(a) \\
& = & \zeta(x,y,z)\inv\e(y, z)(\sigma a) \\
& = & \zeta(x,y,z)\inv\e(y, z)(\sigma\tau a)\inv \\
&=& \zeta_E((x,a), (y, b), (z, c))\inv\sigma\tau \w( (x,a), (y, b),
(z, c))\inv\\
&=& \w'((x,a), (y, b), (z, c))\inv\,.
\end{eqnarray*}
Thus, the map $x \rtimes a \mapsto (x,a)$ defines an isomorphism of
coquasi-bialgebras from $B$ to $\QH{\w'}{E}$.

     \medskip
Now we have equivalences
$$
\C{D^{\zeta_G}(G)}\cong \mathcal{Z}(\Co{B}) \cong
\mathcal{Z}(\Co{\QH{\w'}{E}}) \cong \C{D^{\sigma\tau\w\zeta_E}(E)}
$$
of braided tensor categories. Since $\sigma$ is arbitrary,
$\sigma\tau\w$ can be any element of $\Omega$. This completes the
proof of the Theorem.
\end{proof}

\begin{rem}
   If $A$ is a trivial $H$-module, one can verify directly via rather
extensive computation that
\begin{equation} \label{eq:F}
F = \sum_{\bsa p,q \in H \\ \mu, \nu \in \wh{A} \esa} e(p,\mu) \ot
(1,\nu^{-1}) \ot e(q,\nu) \ot (1,1),
    \end{equation}
    is a gauge transformation in $D^{\zeta_G}(G) \ot D^{\zeta_G}(G)$, and that
    $\varphi\colon D^{\w\zeta_E}(E) \to D^{\zeta_G}(G)_{F, u}$, given by
   \begin{equation}\label{eq:iso}
\vp: e(h,a) \ot (k,b) \mapsto \frac{\epsilon(h,k)(b)}{\epsilon(k,
k^{-1}hk)(b)} \cdot \frac{1}{|A|}\sum_{\chi,\psi \in \wh{A}}
\frac{\chi(b)}{\psi(a)}e(h,\chi) \ot (k,\psi)
    \end{equation}
     with unit
     \begin{eqnarray}
u =  \sum_{\bsa p \in H \\ \mu \in \wh{A} \esa} e(p,\mu) \ot
(1,\mu), \label{eq:u}
\end{eqnarray}
   is an isomorphism of quasi-triangular quasi-bialgebras. We also note
that the special case of Theorem \ref{thmfamily} in which $\zeta$ is
trivial can be deduced without difficulty from some results of
Natale \cite{N}.
\end{rem}

\section{ An Isomorphism of Bialgebras}\label{sectionBI}

\hspace{0.475cm} We next consider a variant of Theorem
\ref{thmfamily}, which however is more limited in scope. We take $G$
to be a group with $N \unlhd G$ a normal subgroup of \emph{index
$2$}. Denote projection onto the quotient as
\begin{eqnarray*}
G \longrightarrow G/N, \ \ x \mapsto \overline{x}.
\end{eqnarray*}

\hspace{0.475cm} Let
    $\eta \in Z^3(G/N, \mathbb{C}^*)$ be the $3$-cocycle which represents the
    nontrivial cohomology class of $G/N = \BZ_2$, so that
    \begin{eqnarray*}
\eta(\overline{x}, \overline{y}, \overline{z}) = \left \{
\begin{array}{ll}
    -1 &   \mbox{if} \  \overline{x}, \overline{y},  \overline{z} \neq 1,  \\
    1 & \mbox{otherwise}.
    \end{array}
    \right.
\end{eqnarray*}
$\eta$ is certainly an abelian $3$-cocycle, that is the associated
$2$-cocycles $\hat{\theta}_{\overline{g}}$ are \emph{coboundaries}.
Thus
\begin{eqnarray*}
\hat{\theta}_{\overline{g}} (\overline{x}, \overline{y}) =
f_{\overline{g}} (\overline{x}) f_{\overline{g}} (\overline{y})  /
f_{\overline{g}} (\overline{x} \overline{y}) =  \left \{
\begin{array}{ll}
    -1 &   \mbox{if} \  \overline{g}, \overline{x},  \overline{y} \neq 1,  \\
    1 & \mbox{otherwise},
    \end{array}
    \right.
\end{eqnarray*}
    where $f_{\overline{g}}$ satisfies
    \begin{eqnarray*}
f_{\overline{g}}(\overline{x}) = \left \{ \begin{array}{ll}
    i &   \mbox{if} \  \overline{g}, \overline{x} \neq 1,  \\
    1 & \mbox{otherwise}.
    \end{array}
    \right.
\end{eqnarray*}

\begin{thm}\label{bialgiso}
     Let the notation be as above, with $\eta' =\infl_{G/N}^G \eta$
     the inflation of $\eta$
to $G$. For any $\omega \in Z^3(G, \mathbb{C}^*)$, the map
    $ \phi: D^{\omega }(G) \rightarrow D^{\omega \eta'}(G)$ given by
$$
\phi: e(g) \otimes x \mapsto f_{\overline{g}}(\overline{x})^{-1}
e(g) \otimes x
$$
is an isomorphism of bialgebras.
    \end{thm}
\begin{proof}
     Let $\theta_g, \theta'_g$ and $\gamma_g, \gamma'_g$ denote the $2$-cochains
    (which are in fact $2$-cocycles) associated with $\omega$ and $\omega \eta'$
    respectively (cf. (\ref{thetaformula}), (\ref{gammaformula})). Clearly
    $$
\theta'_g(x,y) = \theta_g(x, y) \hat{\theta}_{\overline{g}}
(\overline{x}, \overline{y}).
    $$
    Similarly,
     $$
\gamma'_g(x,y) = \gamma_g(x, y) \hat{\theta}_{\overline{g}}
(\overline{x}, \overline{y})^{-1}.
    $$
     Then for any $g, x, h, y \in G$, we have
$$
\begin{aligned}
\phi((e(g) \otimes x)\cdot (e(h) \otimes y)) & = \theta_g(x,y)
\delta_{g^x}^h\phi(e(g) \otimes xy)\\
&=\theta_g(x,y) f_{\overline{g}}(\overline{xy})^{-1} \delta_{g^x}^h
e(g) \otimes xy \\
    &= \theta_g(x,y) \hat{\theta}_{\overline{g}} (\overline{x},
\overline{y})  f_{\overline{g}} (\overline{x})^{-1}
f_{\overline{g}} (\overline{y})^{-1}  \delta_{g^x}^h e(g) \otimes xy \\
&= \theta_g'(x,y)  f_{\overline{g}} (\overline{x})^{-1}
f_{\overline{h}} (\overline{y})^{-1}  \delta_{g^x}^h e(g) \otimes xy \\
&= \phi(e(g) \otimes x)\cdot\phi(e(h) \otimes y).
\end{aligned}
$$
    Similarly,
$$
\begin{aligned}
(\phi \otimes \phi)\Delta(e(g) \otimes x) &
    = \sum_{ab=g}
\gamma_x(a,b) \phi(e(a) \otimes x)
\otimes \phi(e(b)\otimes x) \\
&=\sum_{ab=g} \gamma_x(a,b) f_{\overline{a}}(\overline{x})^{-1}
f_{\overline{b}}(\overline{x})^{-1} e(a) \otimes x
\otimes e(b)\otimes x  \\
    &=\sum_{ab=g}
\gamma_x(a,b) f_{\overline{x}}(\overline{a})^{-1}
f_{\overline{x}}(\overline{b})^{-1} e(a) \otimes x
\otimes e(b)\otimes x  \\
     &=\sum_{ab=g}
\gamma_x(a,b)
f_{\overline{x}}(\overline{g})^{-1}\hat{\theta}_{\overline{x}}
(\overline{a}, \overline{b})^{-1}  e(a) \otimes x
\otimes e(b)\otimes x  \\
     &= \sum_{ab=g}
\gamma_x'(a,b)  f_{\overline{g}}(\overline{x})^{-1} e(a) \otimes x
\otimes e(b) \otimes x \\
&=\Delta(\phi(e(g) \otimes x)). \qedhere
\end{aligned}
$$
\end{proof}
\vspace{.4cm}
    Despite its similarity to Theorem \ref{thmfamily}, we will later
show that Theorem \ref{bialgiso} cannot be improved
    in the sense that  the two twisted doubles that occur
    are generally \emph{not} gauge equivalent.

    \begin{cor}\label{corbialgiso}
    Let notation be as in Theorem \ref{bialgiso}. Then  $D^{\omega }(G)$
and $D^{\omega \eta'}(G)$
    have isomorphic fusion algebras.
\end{cor}

\section{Extra-Special 2-Groups}\label{sectionES}

\hspace{0.475cm} Let $l \geq 1$ be an integer and $V \cong
\BZ_2^{2l}$. Throughout this Section, $Q$ denotes an
\emph{extra-special} $2$-group, defined via a central extension of
groups
\begin{eqnarray}\label{exspecdef}
1 \longrightarrow Z \longrightarrow Q
\stackrel{\pi}{\longrightarrow} V  \longrightarrow 1
\end{eqnarray}
together with the additional constraint that $\BZ_2 \cong Z = Z(Q)$.
Then $Z$ coincides with the commutator subgroup $Q'$ and $|Q| =
2^{2l+1}$. Note that $Q$ has exponent exactly $4$. $l$ is called the
\emph{width} of $Q$, and for each value of $l$, $Q$ is isomorphic to
one of just two groups, denoted by $Q_+, Q_{-}$. If $l=1$ then $Q_+
\cong D_8$ (dihedral group)  and $Q_{-} \cong Q_8$ (quaternion
group). In general $Q_+$ is isomorphic to the central product of $l$
copies of $D_8$ and $Q_{-}$ to the central product of $l-1$ copies
of $D_8$ and one copy of $Q_8$. The squaring map $x \mapsto x^2, x
\in Q,$ induces a nondegenerate quadratic form $q:V \rightarrow
\BF_2,$ and $Q$ has type $+$ or $-$ according as the Witt index $w$
of $q$ is $l$ or $l-1$ respectively.

    \medskip
We need to develop some facts concerning $H^3(Q, \mathbb{C}^*)$. For
background on the cohomology of extra-special groups one may consult
\cite{Q}, \cite{HK} and \cite{BC}.  We use the following notation:
for an abelian group $B$, $\Omega B = \{ b \in B | b^2 = 1 \}$;
     $\Omega B$ is an elementary abelian $2$-group.\\

    \begin{lem}\label{cohomgroups} The following hold:
\begin{eqnarray*}
&&a) H^3(Q, \mathbb{C}^*) = \BZ_2^N \oplus \BZ_4 \ (Q \not\cong
Q_8), \ N = {2l \choose 1}+ {2l \choose 2}+ {2l \choose
3} - 1;  \\
&& b) \   H^3(D_8, \mathbb{C}^*) =
\BZ_2^2 \oplus \BZ_4; \\
&&c) \  H^3(Q_8, \mathbb{C}^*) =  \BZ_8;
\end{eqnarray*}
     \end{lem}
\begin{proof}
Part a) is proved in \cite{HK}, and b) is a special case. It is
well-known (\cite{CE}, Theorem $11.6$) that $Q_8$ has periodic
cohomology. In particular $H^4(Q_8, \BZ)$ is cyclic of order $8$,
and c) holds.
\end{proof}

\begin{lem}\label{lemmainflation} Suppose that $Q \not \cong D_8.$
    Then
\begin{eqnarray*}
\Omega H^3(Q, \mathbb{C}^*) =  \infl_V^Q H^3(V, \mathbb{C}^*).
\end{eqnarray*}
\end{lem}
\begin{proof}
Let $\beta$ and $\beta'$  be the connecting maps (Bocksteins)
associated with the short exact   sequences of coefficients
\begin{eqnarray*}
&& 0 \longrightarrow \BZ \stackrel{2}{\rightarrow} \BZ
\longrightarrow
\BF_2 \longrightarrow 0, \\
&& 0 \longrightarrow \BZ {\longrightarrow} \mathbb{R}
\stackrel{exp}{\longrightarrow} \mathbb{S}^1 \longrightarrow 1.
\end{eqnarray*}
   From the associated long exact sequences in cohomology we obtain a
commuting diagram
\begin{eqnarray}\label{bockdiag}
\begin{array}{ccccc}
    H^3(Q, \BF_2) & \stackrel{\beta_Q}{\longrightarrow}  \Omega
H^4(Q, \BZ) &
    \stackrel{\beta'}{\longleftarrow} & \Omega H^3(Q, \mathbb{C}^*) \\
     \uparrow &      \uparrow   && \uparrow   \\
H^3(V, \BF_2) & \stackrel{\beta_V}{\longrightarrow}  H^4(V, \BZ) &
\stackrel{\beta'}{ \longleftarrow} &   H^3(V, \mathbb{C}^*)
    \end{array}
\end{eqnarray}
where vertical maps are inflations, $\beta'$ is an isomorphism, and
$\beta$ a surjection. It therefore suffices to show that the
leftmost vertical map is surjective. In \cite{Q}, Quillen describes
the cohomology ring $H^*(Q, \BF_2)$
    as a tensor product
    \begin{eqnarray*}
    H^*(Q, \BF_2) = S(V)/J \otimes \BF_2[\zeta]
\end{eqnarray*}
where the left factor $S(V)/J$ (a certain quotient of the cohomology
$H^*(V, \BF_2) = S(V)$) coincides with $\infl_V^Q H^*(V, \BF_2)$,
and $\zeta$ is a cohomology class of degree $2^{2l-w}$. From this we
see that if the leftmost vertical map in (\ref{bockdiag}) is
\emph{not} surjective then $\BF_2[\zeta]$ has a nonzero element of
degree 2. Then $2l-w \leq 1$, and this can only happen if $l = w =
1$, i.e. $Q \cong D_8$. This completes the proof of the Lemma.
\end{proof}

There is a distinguished element $\nu \in H^3(V, \BF_2)$ that we
will need. Namely, pick a basis of $V$ and let $\lambda_1, ...,
\lambda_{2l} \in \mbox{Hom}(V, \BF_2)$ be the dual basis. Let the
quadratic form $q$ corresponding to $Q$ be $q = \sum_{i, j}
a_{ij}x_ix_j$, and set
\begin{eqnarray}\label{nudef}
\nu = \sum_{i,j} a_{ij} \lambda_i^2 \lambda_j \in S^3(V) = H^3(V,
\BF_2).
\end{eqnarray}

    \begin{lem}\label{lemmaradelt}  Suppose that $Q \not \cong D_8$.
    Then
    $H^4(Q, \BZ) \ \cap$ rad $H(Q, \BZ)$
     is cyclic, and the unique element of order $2$ that it contains is
$\infl_V^Q \beta_V(\nu) = \beta_Q(\infl_V^Q  \nu)$.
     \end{lem}
     \begin{proof}
      From Lemma \ref{lemmainflation} we know that any element
      $\zeta \in \Omega H^4(Q, \BZ) \ \cap$ rad $H(Q,
\BZ)$  satisfies
      \begin{eqnarray*}
      &&\zeta =  \beta_Q( \infl_V^Q \delta) =   \infl_V^Q
\beta_V(\delta), \\
      &&\mbox{Res}_C^Q  \zeta = 0, \forall \ C \subseteq Q \ \mbox{of
order $2$},
            \end{eqnarray*}
             for some
       $\delta \in H^3(V, \BF_2)$.

     \medskip
     We claim that $\Res_D^V(\delta)=0$
for every subgroup $D \subseteq V$ of order 2 generated by a {\em
singular vector} $x$ (i.e. $q(x) = 0$). First note that with the
notation of (15), the singularity of
     $x$ means that $P=\pi^{-1}D$ is  the direct
product of $Z$ and $C$ for some subgroup $C \subseteq Q$
     of order $2$.
Consider the commutative diagram
$$
\xymatrix{
0 \ar[r]& H^3(C, \BF_2)\ar[r]^{\beta_C}& H^4(C, \BZ)\\
    &  H^3(P, \BF_2) \ar[u]^{\Res}
\ar[r]^{\beta_P} & H^4(P, \BZ)
\ar[u]^{\Res}\\
H^3(P/Z, \BF_2)\ar[ur]^{\infl}  &  H^3(Q, \BF_2) \ar[u]^{\Res}
\ar[r]^{\beta_Q} & H^4(Q, \BZ)
\ar[u]^{\Res}\\
& H^3(V, \BF_2)\ar[lu]^{\Res}\ar[u]^{\infl}\ar[r]^{\beta_V} & H^4(V,
\BZ)\ar[u]^{\infl}}
$$
Since $\Res_C^Q \circ \beta_Q \circ \infl_V^Q(\delta) = \Res_C^Q
\zeta =0$, we find $\Res_C^Q \circ \infl_V^Q(\delta)=0$ by the
injectivity of $\beta_C$ and the commutative diagram. Notice that
the composition $\Res_{C}^P \circ \infl_D^P$ is an isomorphism. It
follows from the commutative diagram that $\Res_{D}^V(\delta)=0$.

    \medskip
     The classes $\delta$ with this property are spanned by
     the elements $\lambda_i^2\lambda_j + \lambda_i \lambda_j^2$
     together with $\nu$.  However,   the Bockstein annihilates all of
these elements
     with the exception of $\nu$.
    This proves that $ \Omega H^4(Q, \BZ) \ \cap$ rad $H(Q, \BZ)$ is
    spanned by $\beta_Q(\infl_V^Q  \nu)$, and in particular
    $H^4(Q, \BZ) \ \cap$ rad $H(Q, \BZ)$ is cyclic.
    \end{proof}

    \begin{lem}\label{lemmaD8} Let $\zeta \in H^3(D_8, \mathbb{C}^*)$
    have order $4$, and set $I = \infl_V^{D_8}H^3(V, \mathbb{C}^*)$. Then
    \begin{eqnarray*}
    H^3(D_8, \mathbb{C}^*) = \langle \zeta \rangle \times  I.
    \end{eqnarray*}
    \end{lem}
     \begin{proof}
     In the case $Q = D_8$, the proof of Lemma \ref{lemmainflation}
shows that $|I| = 4$, so after Lemma \ref{cohomgroups}b) it suffices
     to show that $I$ contains no nonidentity squares. But it is easily
checked (cf. \cite{M2})
     that the three nonzero elements of $I$ are inflated from elements
     $\delta \in H^3(V, \mathbb{C}^*)$ with the property that $\delta$
     restricts nontrivially to exactly two subgroups of $V$ of order $2$. Then
     $\infl_V^{D_8} \delta$ must restrict nontrivially to some order $2$
     subgroup of $D_8$, hence cannot be a radical element. This
completes the proof of the Lemma.
     \end{proof}

\begin{lem}\label{cohomrestrict} Let $C \subseteq Q$ be
    cyclic of order $4$, and
    $R =\Res_C^Q H^3(Q, \mathbb{C}^*).$ Then
\begin{eqnarray*}
    |R| = \left \{ \begin{array}{ll}
    2 &  \  \mbox{if $l \geq2$}, \\
     4 &  \  \mbox{if $l = 1$}.
     \end{array} \right.
      \end{eqnarray*}
\end{lem}
\begin{proof}
   First assume that $Q \cong Q_8$. Then a generator of $H^3(Q,
\mathbb{C}^*)$ (cf. Lemma \ref{cohomgroups}a)) restricts to a
generator of $H^3(C, \mathbb{C}^*)$ for any subgroup $C \subseteq Q$
(\cite{CE}), and in particular $|R|=4$ in this case.

Next assume that $l \geq 2$. There is a unique nontrivial square in
$H^3(Q, \mathbb{C}^*),$
    call it $\eta$. Then $\eta$ is described in Lemma \ref{lemmaradelt},
and that result shows that for any subgroup $Q_1 \subseteq Q$
satisfying $Q_1 \cong Q_8$, the element  Res$_{Q_1}^Q \eta$ is
\emph{nontrivial}. Now let $\zeta \in H^3(Q, \mathbb{C}^*)$ have
order $4$. Then $\zeta^2 = \eta$, and it follows from what we have
said that Res$_{Q_1}^Q \zeta$ continues to have order $4$. From the
first paragraph of the proof  it  follows that Res$_C^Q \zeta$ has
order $2$ whenever $C \subseteq Q_1$ is cyclic of order $4$. On the
other hand, since $Q \not \cong D_8$, any cyclic subgroup of $Q$ of
order $4$ is contained in some $Q_1$. This completes the proof of
the Lemma in the case $l \geq 2.$

\medskip
It remains to deal with the case $Q \cong D_8$. We use an argument
based on the Lyndon-Hochschild-Serre spectral sequence associated to
the short exact sequence
\begin{eqnarray*}
1 \rightarrow \BZ_4 \rightarrow D_8 \rightarrow \BZ_2 \rightarrow 1.
\end{eqnarray*}
The $E_2$-term is
\begin{eqnarray*}
E_2^{pq} = H^p(\BZ_2, H^q(\BZ_4, \BZ)).
\end{eqnarray*}
Now
\begin{eqnarray*}
&&E_2^{31} = E_2^{13} = 0, \\
&&E_2^{04} = H^4(\BZ_4, \BZ)^{\BZ_2} =
H^4(\BZ_4, \BZ) = \BZ_4, \\
&&E_2^{22} = H^2(\BZ_2, H^2(\BZ_4, \BZ)) =
\BZ_2, \\
&&E_2^{40} = H^4(\BZ_2, \BZ) = \BZ_2.
\end{eqnarray*}
Using Lemma \ref{cohomgroups}b) we see that $|H^4(D_8, \BZ)| =
|E_2^{04}||E_2^{22}||E_2^{40}| = 16$, whence $E_2^{04} =
E_{\infty}^{04}, E_2^{22} = E_{\infty}^{22}, E_2^{40}
=E_{\infty}^{40}.$ In particular
\begin{eqnarray*}
\mbox{Im(Res:} \ H^4(D_8, \BZ) \rightarrow H^4(\BZ_4, \BZ)) =
E_{\infty}^{04} = H^4(\BZ_4, \BZ).
\end{eqnarray*}
This completes the proof of  the Lemma.
\end{proof}

We now apply these results to  twisted quantum doubles. To begin
with, by combining Lemma \ref{lemmainflation} with Theorem
\ref{thmfamily} we immediately obtain

\begin{thm}\label{exspecqdbleisom} Let $Q$ be an extra-special $2$-group
of width $l$ with $Q \not \cong D_8$, and let $E$ be an elementary
abelian $2$-group of order $2^{2l+1}$. For any $\eta \in \Omega
H^3(Q, \mathbb{C}^*)$, there are $\mu \in H^3(E, \mathbb{C}^*)$ and
an
     equivalence of braided tensor categories
\begin{eqnarray}\label{gequiv}
\C{D^{\mu}(E)} \cong \C{D^{\eta}(Q)}.
\end{eqnarray}
\end{thm}

    After the last Theorem the question arises: what can one say about
$D^{\eta}(Q)$ in the case where $\eta$ has order \emph{greater} than
$2$? We will show that an  equivalence of the type (\ref{gequiv}) is
not possible in this case. For this we need to compute some gauge
invariants of the twisted quantum doubles. We use the
\emph{Frobenius-Schur exponent}, which is particularly accessible.

\medskip
    First recall
(\cite{MN2},  \cite{NS1}, \cite{NS2}) the Frobenius-Schur index
$\nu(\chi)$ and their higher order analogs $\nu^{(n)}(\chi)$
associated to an irreducible character $\chi$ of a quasi-Hopf
algebra and a positive integer $n$. These are indeed gauge
invariants (loc. cit.) The smallest positive integer $f$ such that
$\nu^{(f)}_\chi = \chi(1)$ for all irreducible characters $\chi$ of
$D^{\omega}(G)$ is called the the {\em Frobenius-Schur exponent} of
$D^{\omega}(G)$. The existence of $f$ for general quasi-Hopf
algebras is shown in \cite{NS3}, and it is therefore also a gauge
invariant. The following explicit formula for the Frobenius-Schur
exponent $f = f(G, \omega)$ of a twisted quantum double
$D^{\omega}(G)$ is known (loc. cit.):
\begin{eqnarray}\label{FSexp}
f = \mbox{lcm} \left( \  |C|  \ |\omega_C| \  \right).
\end{eqnarray}
The notation is as follows: $C$ ranges over cyclic subgroups of $G$,
$\omega_C$ is the restriction of $\omega$ to $C$, and $|\omega_C|$
is the order of the corresponding cohomology class. We note that in
\cite{Nat05},   Natale defined the {\em modified exponent }$\exp_\w
G$ of a group $G$ endowed with a 3-cocycle $\w$ by same formula
(\ref{FSexp}),
    and proved that $\exp_\w G$ is a gauge
invariant of $D^\w(G)$ for groups of odd order. In \cite{NS3} it is
shown that generally that $f$ is a gauge invariant of $D^\w(G)$ and
that $f=\exp(D^\w(G))$ or  $2\exp(D^\w(G))$.

We now have

\begin{thm}\label{FSexpcomp} The following hold for \emph{any}
$3$-cocycle $\omega$ on one of the indicated groups  ($Q$ is an
extra-special group of width $l$):
    \begin{eqnarray*}
    && a) \ f(\BZ_2^l, \omega) = \left \{ \begin{array}{ll}
4 & \mbox{if} \ \omega \not = 1, \\
2 & \mbox{if $ \omega = 1$,} \ l \geq 1;
\end{array} \right.\\
&& b) \ f(Q, \omega) = \left \{ \begin{array}{ll}
4 & \mbox{if} \ \omega^2 = 1, \  l \geq 2, \\
8 & \mbox{if $\omega^2 \not = 1, \  l  \geq 2$};
\end{array} \right.\\
    && c) \ f(Q_8, \omega) = \left \{ \begin{array}{lll}
4 & \mbox{if} \ \omega^2 = 1, \\
8 & \mbox{if} \ \omega \ \mbox{has order $4$}, \\
16 & \mbox{if} \ \omega \ \mbox{has order $8$}; \\
\end{array}
\right.  \\
    && d) \ f(D_8, \omega) = \left \{ \begin{array}{lll}
4 & \mbox{if} \ \omega \in  \infl_V^{D_8}H^3(V, \mathbb{C}^*), \\
8 & \mbox{if} \ \omega \in \Omega H^3(D_8, \mathbb{C}^*) \setminus
\infl_V^{D_8}H^3(V, \mathbb{C}^*), \\
16 & \mbox{if} \ \omega \ \mbox{has order $4$}. \\
\end{array}
\right.
    \end{eqnarray*}
\end{thm}
\begin{proof}

If $\omega$ is trivial then by (\ref{FSexp}) $f$ coincides with the
exponent of the group in question. Suppose that $\omega$ has order
$2$ in the elementary abelian $2$-group case. Then from \cite{M2}
that there is a subgroup $C$ of order $2$ such that $\omega_C$ is
nontrivial.
    Part a) now follows.
    Parts b) and c) follow immediately from Lemma \ref{cohomrestrict}.

    \medskip
    It remains to handle the case $Q = D_8$. Let $C$ be the cyclic
subgroup of order $4$ and let $I =\infl_V^{D_8}H^3(V,
\mathbb{C}^*)$.   By Lemma \ref{lemmaD8} we have
    $H^3(D_8, \mathbb{C}^*) = I\langle \zeta \rangle$ for any $\zeta$ of
order $4$.
    If $\omega \in I$ then
    $\omega_C$ is necessarily trivial, and therefore
    $f(D_8, \omega) = 4$ for such $\omega$. If $\omega$ has order $4$ then
    $\omega_C$ also has order $4$ (Lemma \ref{cohomrestrict}), so $f =
16$ in this case.
    Finally, $\zeta^2_C$ has order $2$ for $\zeta$ of order $4$. Since
$\zeta^2 \notin I$ we conclude that $\omega_C$ has order $2$
    whenever $\omega \in \Omega H^3(D_8, \mathbb{C}^*) \setminus I$, so
that $f = 8$ in this case. This completes the proof of all parts
    of the Theorem.
\end{proof}

\section{Twisted Quantum Doubles of Dimension 64}\label{section64}

\hspace{0.475cm} We have seen in the previous Section that the
properties of the groups $H^3(Q, \mathbb{C}^*)$ for $Q = Q_8, D_8$
are exceptional in several ways. Thus the same is true for the
corresponding twisted quantum doubles.
    In this Section we consider the problem of
understanding the gauge equivalence classes of the quantum doubles
$D^{\omega}(G)$ where $G$ has order $8$. In \cite{MN1},  the same
question was treated for \emph{abelian} groups $G$ and
\emph{abelian} cocycles, that is cocycles $\omega$ for which
$D^{\omega}(G)$ is commutative. Here we consider the case when
$D^{\omega}(G)$ is \emph{noncommutative}. This precludes the two
groups $\BZ_8$ and $\BZ_2 \times \BZ_4$, leaving the groups $E_8 =
\BZ_2^3, Q_8$ and $D_8$ to be considered. In this Section we denote
these three groups by $E, Q$ and $D$ respectively. Now $E$ has $64$
nonabelian $3$-cohomology classes, while $D$ and $Q$ have $16$ and
$8$ classes respectively (cf. Remark \ref{cohomgroups}).
    So our task is to sort the
$88$ resulting twisted quantum doubles into gauge equivalence
classes.  We will see that there are at least 8, and no more than
20, gauge equivalence classes.

\vspace{0.4cm} A significant reduction in the problem is achieved by
considering automorphisms of the three groups in question. That is
because it is easy to see (\cite{MN1}, Remark 2.1(iii)) that
automorphisms preserve gauge equivalence classes. Consider first the
group $E$. One knows that $H = H^3(E, \mathbb{C}^*)$ has order $2^7$
and exponent $2$. Of the $128$ cohomology classes in $H$, the
abelian classes form a subgroup of order $64$ generated by the Chern
classes of characters of $E$ (\cite{MN1}, Proposition 7.5).

\vspace{0.4cm} The automorphism group of $E$ is the simple group
$SL(3, 2)$ of order $168$, and we regard $H$ as a $7$-dimensional
$\BF_2SL(3, 2)$-module.
     We can understand the structure of $H$
    as follows (for more details, see \cite{M2}).
    A cohomology class $\omega \in H$
is characterized by the subset of elements $1 \neq g \in E$ for
which the restriction Res$_{\langle g \rangle}^{E}\omega$ is
\emph{nontrivial}. Let us call this set of elements the
\emph{support} of $\omega$, denoted
    supp $\omega$.
     The \emph{weight} of $\omega$ is the cardinality $|\mbox{supp} \
\omega|$ of its support.
     There is an isomorphism of $SL(3, 2)$-modules
     \begin{eqnarray}\label{cohomisom}
H \stackrel{\cong}{\longrightarrow} \BF_2^7, \ \ \omega \mapsto
\mbox{supp} \ \omega.
\end{eqnarray}
So $H$ is a permutation module for $SL(3, 2)$ corresponding to  the
permutation action of $SL(3, 2)$ on the nonzero elements of $E$.
      Thus there is a decomposition
\begin{eqnarray}\label{cohomdecomp}
H = 3 \oplus \overline{3} \oplus 1
\end{eqnarray}
into simple modules. The abelian cocycles are those in the unique
submodule of codimension 1, and they may therefore be alternately
characterized as those classes
    of  \emph{even} weight.
As for the nonabelian classes, the possible weights are $1, 3, 5,
7$, and the number of cohomology classes of each type is $7,35, 21
,1$ respectively. Those of weight $1, 5$ or $7$ form  single $SL(3,
2)$-orbits of size $7,21, 1$ respectively.
    Those of weight $3$ split into two orbits
according as supp $\omega$ is a set of linearly dependent ($7$ of
them) or linearly independent ($28$ of them) elements
    of  $E$. We will utilize this information to  label the
    cohomology classes of $E$.

    \medskip
     Based on what we have said so far, there are (at most) 5 gauge
equivalence classes of
    noncommutative twisted quantum doubles $D^{\omega}(E_8)$, with
representative
    cocycles $\omega_1, \omega_{3i},  \omega_{3d},  \omega_{5},
\omega_7$. Here, the numerical subscript is the
    weight of the cocycle, and the additional subscript $i$ or $d$ at
weight 3 indicates whether
supp $\omega$ is a linearly independent or dependent set respectively. \\

     Next we consider the quantum doubles $D(D)$ and $D(Q)$.
    Let $\epsilon_Q, \epsilon_D$ denote elements
     of $H^2(V, \BZ_2)$ defining $Q$ and $D$ respectively
     as central extensions (\ref{exspecdef}). If $h_1, h_2$ are generators
     of $V$ and $t$ a generator of $\BZ_2$ then we may take both
$\epsilon_Q$
     and $\epsilon_D$ to be (multiplicatively)  bilinear, and
     \begin{eqnarray}
&\epsilon_D(h_1, h_1) = \epsilon_D(h_1, h_2) =  \epsilon_D(h_2, h_2)
= 1,  \epsilon_D(h_2, h_1) = t;
\label{dihedralcocycle} \\
&\epsilon_Q(h_1, h_1) = \epsilon_Q(h_2, h_1) =  \epsilon_Q(h_2, h_2)
= t,  \epsilon_Q(h_1, h_2) = 1. \label{quaternioncocycle}
\end{eqnarray}

    There is a useful isomorphism analogous to
(\ref{cohomisom}) (cf. \cite{M2}), namely
     \begin{eqnarray}\label{autS3iso}
H^3(V, \mathbb{C}^*) \longrightarrow \BF_2^3, \ \ \zeta \mapsto
\mbox{supp} \ \zeta.
\end{eqnarray}
Identifying $E$ with $\langle h_1, h_2 \rangle \times \langle t
\rangle$, we find
    that the $3$-cocycles $\omega_D, \omega_Q \in H^3(E, \mathbb{C}^*)$
associated to $\epsilon_D, \epsilon_Q$ respectively by
(\ref{cocycledef}) satisfy
     \begin{eqnarray*}
\mbox{supp} \ \omega_D &=& \{(h_1h_2, t)\}, \\
\mbox{supp} \ \omega_Q &=& \{(h_1, t), (h_2, t), (h_1h_2, t) \}.
\end{eqnarray*}
By the case $\zeta = 1$ of Theorem \ref{thmfamily} we conclude that
\begin{eqnarray}
D(D) &\sim& D^{\omega_1}(E), \label{firstDequiv}\\
D(Q) &\sim& D^{\omega_{3 i}}(E). \label{firstQequiv}
\end{eqnarray}

With the notation of Lemma \ref{lemmaD8},  that result shows that
$I$ is generated
    by two classes $\alpha_i, i = 1,2$ such that
    \begin{eqnarray*}
\mbox{supp} \ \alpha_i = \{h_i\}.
\end{eqnarray*}
Now set $\alpha_i' = \infl_{E_4}^{E_8} \alpha_i$, and note that
\begin{eqnarray*}
\mbox{supp} \ \alpha_1'\omega_D &=& \{ (h_1, 1), (h_1, t), (h_1h_2, t)  \}, \\
\mbox{supp} \ \alpha_2'\omega_D &=& \{ (h_2, 1), (h_2, t), (h_1h_2, t)   \}, \\
\mbox{supp} \ \alpha_1' \alpha_2' \omega_D &=& \{ (h_1, 1), (h_1,
t), (h_2, 1), (h_2, t), (h_1h_2, t)   \}.
\end{eqnarray*}
According to Theorem \ref{thmfamily} we can conclude that
\begin{eqnarray}
D^{\alpha_1}(D) &\sim& D^{\alpha_2}(D) \sim D^{\omega_{3 i}}(E),
\label{secondDequiv} \\
D^{\alpha_1\alpha_2}(D) &\sim& D^{\omega_5}(E). \label{thirdDequiv}
\end{eqnarray}
Note that the gauge equivalence $D^{\alpha_1}(D) \sim
D^{\alpha_2}(D)$ also follows from the observation that an
involutorial automorphism $a$ of $V$ lifts to an automorphism of
$D$, also denoted as $a$, which commutes with inflation. Then $a$
exchanges $\alpha_1$ and $\alpha_2$ and therefore induces a gauge
equivalence between the corresponding twisted quantum doubles.
Indeed, it can be verified that there is an $\langle a
\rangle$-equivariant decomposition
\begin{eqnarray*}
H^3(D, \mathbb{C}^*) = \langle \alpha_1 \rangle \oplus \langle
\alpha_2 \rangle \oplus \langle \alpha_3 \rangle
\end{eqnarray*}
where $\alpha_3$ is a class of order $4$ which is $a$-invariant (cf.
Lemma \ref{lemmaD8}). As a result, we get the following additional
gauge equivalences induced by the action of $a$ on cohomology:
\begin{eqnarray}
D^{\alpha_1\alpha_3 }(D_8) &\sim& D^{\alpha_2 \alpha_3}(D_8), \nonumber \\
D^{\alpha_1\alpha_3^2}(D_8) &\sim& D^{\alpha_2\alpha_3^2}(D_8),
\label{Dequivs} \\
D^{\alpha_1\alpha_3^3}(D_8) &\sim& D^{\alpha_2\alpha_3^3}(D_8).
\nonumber
    \end{eqnarray}

    The automorphism group of $Q$ acts \emph{trivially} on the relevant
cohomology group, so that no new gauge equivalences can be realized
from automorphisms of $Q$.
    Let $\gamma$ be a generator of $H^3(Q, \mathbb{C}^*)$ (cf. Lemma
    \ref{cohomgroups}c)).  By Lemma \ref{lemmainflation}
    we have  $\infl_{E_4}^{Q}H^3(E_4, \mathbb{C}^*) = \langle \gamma^4 \rangle$.
     We
then obtain, after a calculation,  the following gauge equivalence
that arises from application of Theorem \ref{thmfamily}:
\begin{eqnarray}\label{Qequivs}
D^{\gamma^4}(Q_8) \sim D^{\omega_{3d}}(E_8).
\end{eqnarray}

We have now established
\begin{thm} The 88 noncommutative twisted quantum doubles of $E, D$ and $Q$
fall into  \emph{at most} 20 gauge equivalence classes, namely
(\ref{firstDequiv}),  (\ref{firstQequiv}),    (\ref{thirdDequiv}),
(\ref{Dequivs})   and (\ref{Qequivs}) together
    with the class of $D^{\omega_7}(E_8)$, the 6 classes $D^{\alpha}(D_8)$
    with $\alpha \ne 1$ or $\a_1\a_2$ in $\langle \alpha_1\alpha_2,
\alpha_3 \rangle$,
    and the classes $D^{\gamma^m}(Q), m = 1, 2, 3, 5, 6, 7$.
    \end{thm}

    It is possible that there are
    \emph{less} than
    20 distinct gauge equivalence classes. We will see that there
\emph{at least} 8 such classes.
    To describe this, consider the following sets of cohomology classes:
    \begin{eqnarray*}
&&E:  \left \{ \begin{array}{ll}
    \mu_1: \omega_1, \omega_7, \\
    \mu_2: \omega_{3i},  \omega_{3d}, \omega_5
    \end{array} \right. \\
    &&D: \left \{ \begin{array}{llllll}
     \eta_0: 1,\\
     \eta_1: \alpha_1, \alpha_2, \alpha_1\alpha_2,\\
     \eta_2 :  \alpha_1 \alpha_3, \alpha_2 \alpha_3, \alpha_1
\alpha_3^3,  \alpha_2 \alpha_3^3, \\
     \eta_3 :  \alpha_3, \alpha_3^3, \alpha_1 \alpha_2 \alpha_3,
\alpha_1 \alpha_2 \alpha_3^3,\\
     \eta_4 :   \alpha_1 \alpha_3^2, \alpha_2 \alpha_3^2 ,  \alpha_1
\alpha_2 \alpha_3^2,\\
     \eta_5: \alpha_3^2.
     \end{array} \right. \\
     &&Q: \left \{ \begin{array}{lll}
      \gamma_0 :  1, \gamma^4, \\
     \gamma_1:  \gamma^2, \gamma^6, \\
     \gamma_2 :  \gamma , \gamma^3, \gamma^5, \gamma^7.
     \end{array} \right.
     \end{eqnarray*}
     Here, we have partitioned representatives for the orbits of
$\Aut G \  (G = E, D, Q),$ acting on the relevant $3$-cohomology
classes of $G$, into
    certain subsets.
    \begin{thm} If $G$ is one of the groups $E, D, Q$, a pair of twisted
quantum doubles $D^{\alpha}(G), D^{\beta}(G)$ have the \emph{same}
sets of (higher) Frobenius-Schur indicators if, and only if,
$\alpha$ and $\beta$ both lie in one of the sets $\mu_i, \eta_j,
\gamma_k$. Between them, there are just $8$ distinct sets of higher
Frobenius-Schur indictators.
\end{thm}
    \begin{proof}
    The complete sets of indicators are given in Appendix below.
The Theorem follows
    from this data.
    \end{proof}

    \medskip
      We note in particular the following table of Frobenius-Schur
exponents. The entries of this table follow
      from Theorem \ref{FSexpcomp}.

     $$
\begin{array}{|c|c|}
\hline \mbox{Twisted Quantum Doubles} &  \mbox{Frobenius-Schur
exponents}\\\hline
     D^{\mu_1}(E_8), D^{\eta_2}(E_8), D^{\eta_0}(D_8), D^{\eta_1}(D_8),
D^{\gamma_0}(Q_8)  & 4 \\\hline
     D^{\eta_4}(D_8), D^{\eta_5}(D_8),  D^{\gamma_1}(Q_8) & 8 \\\hline
     D^{\eta_2}(D_8), D^{\eta_3}(D_8), D^{\gamma_2}(Q_8) & 16 \\\hline
\end{array}
$$

\medskip
    Finally, we briefly consider some bialgebra isomorphisms. Let
$\sigma \in H = H^3(E, \mathbb{C}^*)$
    be a cohomology class inflated, as in Theorem 4.1, from the nontrivial class
    of $E/F$ where $F$ is a subgroup of $E$ of index $2$. Hence, supp
$\sigma$ consists of the
    elements in $E \setminus F$. By Theorem 4.1, there is an isomorphism
    of \emph{bialgebras}
    \begin{eqnarray*}
D^{\omega}(E) \cong D^{\omega \sigma}(E)
\end{eqnarray*}
for all cohomology classes $\omega \in H$.
    Using this, easy arguments lead to  the following bialgebra isomorphisms:
    \begin{eqnarray*}
&& D^{\omega_1}(E) \cong D^{\omega_5}(E) \cong D^{\omega_{3i}}(E), \\
&& D^{\omega_7}(E) \cong D^{\omega_{3d}}(E).
\end{eqnarray*}
As we have seen, $ D^{\omega_7}(E)$ and $D^{\omega_{3d}}(E)$,
    for example,  are not gauge equivalent because they have different
sets of Frobenius-Schur indicators.
    Thus they afford an example of a pair of twisted doubles which
    are isomorphic as bialgebras but not gauge equivalent.

\section{Invariance of Ribbon Structure} \label{ribbon}
Let $(\CC, \ot,  I, \Phi)$ be a left rigid braided monoidal
category, where $I$ is the neutral object and $\Phi$ the
associativity isomorphism. Here, we assume that the neutral object
$I$ is \emph{strict}, i.e. $I \ot V = V \ot I=V$ for $V \in \CC$.
For $V \in \CC$, we use the notation $V\du$ for the left dual of $V$
with the dual basis map $\db_V : I \to V \ot V\du$ and the
evaluation map $\ev_V: V\du \ot V \to I$. Note that $(-)\du$ can be
extended to a contravariant monoidal equivalence of $\CC$ with
$I\du=I$, $\ev_I=\db_V=\id_I$. Thus $(-)\bidu$ is a monoidal
equivalence on $\CC$.  See \cite{K} for more details on right
monoidal category and monoidal functor. We follow the notation and
terminology introduced in \cite{NS1} for the discussion to come.

\medskip
For any braiding $c$ on the left rigid monoidal category $\CC$, the
associated Drinfeld isomorphism $u$ is defined by
\begin{multline*}
   u_V  = \left(V \xrightarrow{\db_{V\du} \ot \id} (V\du \ot V\bidu) \ot V
    \xrightarrow{c\ot \id } (V\bidu \ot V\du) \ot V \xrightarrow{\Phi}\right.\\
    \left. V\bidu \ot (V\du \ot V)
    \xrightarrow{\id \ot \ev}
    V\bidu\right)
\end{multline*}
for $V \in \CC$.  In particular, if $\CC$ is \emph{strict}, the
Drinfeld isomorphism satisfies the equation
$$
u_{V\ot W}= (u_V \ot u_W)c_{V, W}\inv c_{W, V}\inv
$$
for $V, W \in \CC$.

    \medskip
Suppose $j$ is a \emph{pivotal structure} on $\CC$ (not necessarily
strict), i.e. $j: \Id \to (-)\bidu$ is an isomorphism of monoidal
functors. Then $v=u\inv \circ j$ is a {\em twist} of the braided
monoidal category $(\CC, c)$, i.e. a natural automorphism of the
identity functor $\Id$ of $\CC$ such that
$$
v_{V \ot W} = (v_V \ot v_W)c_{W, V}c_{V, W}
$$
for all $V, W \in \CC$. If $\CC$ is a \emph{spherical fusion
category over} $\BC$, then $v$ defines a ribbon structure on $\CC$
(cf. \cite{NS3}). We use the notation $(\CC, c, j)$ to denote a
pivotal braided monoidal category $\CC$ with the braiding $c$ and
the pivotal structure $j$.

\begin{lem}\label{lem6.1}
Let $(\CC, c, j)$, $(\DD, c, j)$ be pivotal braided monoidal
categories. If $(\FF, \xi): \CC \to \DD$ is a braided monoidal
equivalence which preserves the pivotal structures, then
$$
\FF(v_V) = v_{\FF V}
$$
for $V \in \CC$.
\end{lem}
\begin{proof}
    From \cite{NS1} it follows that the duality transformation
$\tilde\xi: \FF(V\du) \to
    \FF(V)\du$ of $(\FF, \xi)$ is determined by either of the
    commutative diagrams
    \begin{equation}\label{eq:dualtran}
    \vcenter{
    \xymatrix@1{\FF(I)\ar[d]^-{\FF(\db)} \ar@{=}[r] & I \ar[dd]^-{\db} \\
    \FF(V \ot V\du)\ar[d]^-{\xi\inv}  & \\
    \FF(V) \ot \FF (V\du) \ar[r]^-{\id\ot \tilde \xi}& \FF(V) \ot   \FF(V)\du
    }
    }
    \quad \text{and}\qquad
     \vcenter{
     \xymatrix{\FF(V\du \ot V)\ar[r]^-{\FF(\ev)} \ar[d]^-{\xi\inv} &
\FF(I) \ar@{=}[dd] \\
               \FF(V\du) \ot \FF(V)\ar[d]^-{\tilde\xi\ot \id}  & \\
    \FF(V)\du \ot \FF (V) \ar[r]^-{\ev}& I
    }
    }\,.
    \end{equation}
We observe that following diagram is commutative: {\tiny
$$
\vcenter{\xymatrix@C-7pt@R+10pt{ \FF V \ar[r]^-{\FF(\db \!\ot\!
\id)} \ar[rd]_-{\FF(\db)\ot  \id}
\ar@/_1pc/{[rddd]+(-10,2)}^-{\db\!\ot\!
\id}\ar@/_2pc/{[rdddd]+(-13,2)}_-{\db \!\ot\!  \id} & \FF((V\du
\!\ot\!  V\bidu) \!\ot\!  V) \ar[d]^-{\xi\inv} \ar[r]^-{\FF(c \ot
\id)} & \FF((V\bidu \!\ot\!  V\du) \!\ot\!  V)\ar[r]^-{\FF(\Phi)}
\ar[d]^-{\xi\inv} &\FF(V\bidu \!\ot\!  (V\du \!\ot\!  V))
\ar[d]^-{\xi\inv} \ar[r]^-{\FF(\id \!\ot\!  \ev)} &
\FF(V\bidu) \ar[ddd]^-{\tilde\xi}\\
& \FF(V\du \!\ot\!  V\bidu) \!\ot\!  \FF V \ar[d]^-{\xi\inv \ot \id}
\ar[r]^-{\FF(c)\ot  \id} & \FF(V\bidu \!\ot\!  V\du)\!\ot\! \FF V
\ar[d]^-{\xi\inv  \ot  \id} & \FF(V\bidu) \!\ot\!  \FF(V\du \!\ot\!
V)\ar[d]^-{\id \ot \xi\inv }
\ar[ru]_-{\id \!\ot  \FF(\ev)}\\
& (\FF(V\du) \!\ot\!  \FF(V\bidu))\!\ot\!  \FF V\ar[d]^-{\id \ot
\tilde\xi  \ot  \id} \ar[r]^-{c \ot  \id} & (\FF(V\bidu) \!\ot\!
\FF(V\du))\!\ot\!  \FF V \ar[d]^-{\tilde\xi  \ot  \id \ot  \id}
\ar[r]^-{\Phi} & \FF(V\bidu) \!\ot\!  (\FF(V\du) \!\ot\!  \FF V)
\ar[d]^-{\tilde\xi
\ot  \id \ot  \id} & \\
& (\FF(V\du) \!\ot\!  \FF(V\du)\du)\!\ot\!  \FF V\ar[r]^-{c\ot  \id}
\ar[d]^-{\tilde\xi \ot  (\tilde \xi\du)\inv \ot  \id } &
(\FF(V\du)\du \!\ot\!  \FF(V\du))\!\ot\!  \FF V\ar[r]^-{\Phi}
\ar[d]^-{(\tilde \xi\du)\inv \ot  \tilde\xi \ot  \id } &
\FF(V\du)\du \!\ot\! (\FF(V\du) \!\ot\!  \FF V)\ar[d]^-{(\tilde
\xi\du)\inv \ot  \tilde\xi \ot  \id }
& \FF(V\du)\du \ar[d]^-{(\tilde\xi\du)\inv} \\
& (\FF(V)\du \!\ot\!  \FF(V)\bidu) \!\ot\!  \FF V \ar[r]^-{c \ot
\id} & (\FF(V)\bidu\!\ot\!  \FF(V)\du)\!\ot\!  \FF V\ar[r]^-{\Phi} &
\FF(V)\bidu\!\ot\!  (\FF(V)\du\!\ot\!  \FF V) \ar[r]^-{\id \!\ot\!
\ev} & \FF(V)\bidu } }.
$$}
Commutativity of the middle rectangles are  consequences of either
the properties
    of  the braided monoidal equivalence $(\FF, \xi)$,
or the naturality of $c$ and $\Phi$. Commutativity of the two
triangles at the upper corners  follows from properties of the
coherence map $\xi$, and that of the lower left triangle follows
from properties of left duality. The
   commutativity of remaining two polygons on both sides follow from
\eqref{eq:dualtran}.

    \medskip
   Note that the top edge is $\FF(u_V)$. The commutativity of the above
diagram implies that
$$
(\tilde \xi\du)\inv \tilde \xi\FF(u_V) = u_{\FF(V)}\,.
$$
Since $(\FF, \xi)$ also preserves pivotal structures (cf.
\cite{NS3}),
$$
(\tilde \xi\du)\inv \tilde \xi = j \FF(j\inv)
$$
and hence
$$
\FF(v_V) = v_{\FF V}\,. \qedhere
$$
\end{proof}

Suppose that $\CC=\Cf{H}$ is the tensor category of
finite-dimensional modules over a semisimple quasi-Hopf algebra $H$.
Then $\CC$ admits a canonical pivotal structure given by a trace
element $g$ of $H$, namely
$$
j_V(x)(f) = f(g\inv x)
$$
for $x \in V$ and $f \in V\du$ (cf. \cite{MN2}, \cite{NS2},
\cite{ENO}). If, in addition, $H$ admits a universal
$\mathcal{R}$-matrix, then $\CC$ is a braided spherical fusion
category. The Drinfeld isomorphism $u_V : V \to V\bidu$ is given by
$$
u_V(x)(f)=f(ux)
$$
where $u$ is the Drinfeld element and the associated ribbon
structure $v_V$ is given by the multiplication of the element $v=(gu)\inv$.\\

\begin{prop}
    Let $H$, $K$ be semisimple, braided, quasi-Hopf algebras. If $(\FF,
\xi): \Cf{H} \to \Cf{K}$ is a braided monoidal equivalence, then
    $$
     \FF(v_V) = v_{\FF(V)}
    $$
    for $V \in  \Cf{H} $.
\end{prop}
\begin{proof}
    By \cite{NS2},  $(\FF, \xi)$ preserves the canonical pivotal
structures. Thus the result follows immediately from Lemma
\ref{lem6.1}.
\end{proof}

Note that $v_V$ is a scalar for any simple $H$-module $V$. Thus if
$(\FF, \xi): \Cf{H} \to \Cf{K}$ is an equivalence of $\BC$-linear
braided monoidal categories, then
$$
v_V = v_{\FF(V)} \quad\text{and}\quad \nu_n(V) = \nu_n(\FF(V))
$$
for all positive integer $n$ and simple objects $V \in \Cf{H}$.

    \medskip
In case $H=D^\w(G)$, it is shown in \cite{MN2} that  the trace
element $g$ of $D^\w(G)$  is given by
$$
g=\sum_{x \in G} \w(x, x\inv, x) e(x) \ot 1.
$$
By \cite{AC92},
$$
   u=\sum_{x\in G} \w(x, x\inv, x)^{-2}
e(x) \ot x\inv,
$$
is the Drinfeld element for the associated braiding of $D^\w(G)$.
Hence
$$
v =(gu)\inv=\sum_{x\in G} \w(x, x\inv, x)\inv e(x) \ot x\inv =
\sum_{x \in G} e(x) \ot x
$$
defines the twist or ribbon structure associated with the underlying
canonical pivotal structure and braiding of $\Cf{D^\w(G)}$. With a
different convention, the formula of the
   ribbon element of $D^\w(G)$ is also shown \cite{AC92}. \\

   We can now use the formula for $v$ to compute the corresponding
    scalar $v_V$ for each simple module $V$ of a 64-dimensional twisted
   double $H$. If $\chi$ is the irreducible character afforded by $V$,
the scalar $v_\chi = \chi(v)/\chi(1)$ is equal to $v_V$.
   The sequence of higher indicators $\nu_\chi^{(n)}$ and $v_\chi$ for
each irreducible character $\chi$ of $H$ are presented in Tables  in
the Appendix. It follows from this data that the 20 classes of
twisted doubles $H$ of dimension 64 identified in Section
\ref{section64} have \emph{distinct} sets of sequences. Since both
higher indicators and $v_\chi$ are preserved by a braided tensor
equivalence, we have

   \begin{thm}
     There are exactly 20 gauge equivalence classes of quasi-triangular
quasi-bialgebras
    among the 64-dimensional
     noncommutative twisted quantum doubles of finite groups. \ \ \ \ \
\ \ \ \ \ \ \ \  \ \ \ \ \   \ \ \ \ \ \ \  $\Box$
   \end{thm}

\section{Appendix}
In the following tables, $m$ denotes the multiplicity of the
sequence, and $\kappa$ is
a primitive 16th complex root of 1 such that $\kappa^4=i$.\\ \\
\begin{tabular}{l}
\bf{Frobenius-Schur Exponent 16} \medskip\\
$
\setlength{\arraycolsep}{1pt}
\renewcommand{\arraystretch}{1}
\begin{array}{|c|c|cccccccccccccc|c|c|}
\hline H &  \chi(1) & \nu^{(2)} & \nu^{(3)}& \nu^{(4)} & \nu^{(5)}&
\nu^{(6)} & \nu^{(7)}& \nu^{(8)} & \nu^{(9)} & \nu^{(10)}&
\nu^{(11)} & \nu^{(12)}& \nu^{(13)} & \nu^{(14)}& \nu^{(15)}& v_\chi
& m \\
\hline
   & 1 &  1 &  1 &  1 &  1 &  1 &  1 &  1 &  1 &  1 &  1 &  1 &  1 &  1
&  1 & 1 & 1 \\
D^{\g}(Q_8) & 1 &  1 &  0 &  1 &  0 &  1 &  0 &  1 &  0 &  1 &  0 &
1 &  0 &  1 &  0 & 1 & 3 \\
   & 1 &  -1 &  0 &  1 &  0 &  -1 &  0 &  1 &  0 &  -1 &  0 &  1 &  0 &
-1 &  0 & -i & 4 \\
   & 2 &  -1 &  0 &  2 &  0 &  -1 &  0 &  2 &  0 &  -1 &  0 &  2 &  0 &
-1 &  0 & 1 & 1 \\
   & 2 &  1 &  0 &  2 &  0 &  1 &  0 &  2 &  0 &  1 &  0 &  2 &  0 &  1
&  0 & i & 1 \\
   & 2 &  -1 &  0 &  1 &  0 &  -1 &  0 &  0 &  0 &  -1 &  0 &  1 &  0 &
-1 &  0 & \ol\kappa^7 & 3 \\
   & 2 &  -1 &  0 &  1 &  0 &  -1 &  0 &  0 &  0 &  -1 &  0 &  1 &  0 &
-1 &  0 & \kappa & 3 \\
   & 2 &  1 &  0 &  1 &  0 &  1 &  0 &  0 &  0 &  1 &  0 &  1 &  0 &  1
&  0 & \ol\kappa^3 & 3 \\
   & 2 &  1 &  0 &  1 &  0 &  1 &  0 &  0 &  0 &  1 &  0 &  1 &  0 &  1
&  0 & \kappa^5 & 3 \\
   \hline
   & 1 &  1 &  1 &  1 &  1 &  1 &  1 &  1 &  1 &  1 &  1 &  1 &  1 &  1
&  1 & 1 & 1 \\
D^{\g^3}(Q_8) & 1 &  1 &  0 &  1 &  0 &  1 &  0 &  1 &  0 &  1 &  0
&
1 &  0 &  1 &  0 & 1 & 3 \\
   & 1 &  -1 &  0 &  1 &  0 &  -1 &  0 &  1 &  0 &  -1 &  0 &  1 &  0 &
-1 &  0 & i & 4 \\
   & 2 &  -1 &  0 &  2 &  0 &  -1 &  0 &  2 &  0 &  -1 &  0 &  2 &  0 &
-1 &  0 & 1 & 1 \\
   & 2 &  1 &  0 &  2 &  0 &  1 &  0 &  2 &  0 &  1 &  0 &  2 &  0 &  1
&  0 & -i & 1 \\
   & 2 &  -1 &  0 &  1 &  0 &  -1 &  0 &  0 &  0 &  -1 &  0 &  1 &  0 &
-1 &  0 & \ol\kappa^5 & 3 \\
   & 2 &  -1 &  0 &  1 &  0 &  -1 &  0 &  0 &  0 &  -1 &  0 &  1 &  0 &
-1 &  0 & \kappa^3 & 3 \\
   & 2 &  1 &  0 &  1 &  0 &  1 &  0 &  0 &  0 &  1 &  0 &  1 &  0 &  1
&  0 & \ol\kappa & 3 \\
   & 2 &  1 &  0 &  1 &  0 &  1 &  0 &  0 &  0 &  1 &  0 &  1 &  0 &  1
&  0 & \kappa^7 & 3 \\
\hline
   & 1 &  1 &  1 &  1 &  1 &  1 &  1 &  1 &  1 &  1 &  1 &  1 &  1 &  1
&  1 & 1 & 1 \\
D^{\g^5}(Q_8) & 1 &  1 &  0 &  1 &  0 &  1 &  0 &  1 &  0 &  1 &  0
&
1 &  0 &  1 &  0 & 1 & 3 \\
   & 1 &  -1 &  0 &  1 &  0 &  -1 &  0 &  1 &  0 &  -1 &  0 &  1 &  0 &
-1 &  0 & -i & 4 \\
   & 2 &  -1 &  0 &  2 &  0 &  -1 &  0 &  2 &  0 &  -1 &  0 &  2 &  0 &
-1 &  0 & 1 & 1 \\
   & 2 &  1 &  0 &  2 &  0 &  1 &  0 &  2 &  0 &  1 &  0 &  2 &  0 &  1
&  0 & i & 1 \\
   & 2 &  -1 &  0 &  1 &  0 &  -1 &  0 &  0 &  0 &  -1 &  0 &  1 &  0 &
-1 &  0 & \ol\kappa^3 & 3 \\
   & 2 &  -1 &  0 &  1 &  0 &  -1 &  0 &  0 &  0 &  -1 &  0 &  1 &  0 &
-1 &  0 & \kappa^5 & 3 \\
   & 2 &  1 &  0 &  1 &  0 &  1 &  0 &  0 &  0 &  1 &  0 &  1 &  0 &  1
&  0 & \ol\kappa^7 & 3 \\
   & 2 &  1 &  0 &  1 &  0 &  1 &  0 &  0 &  0 &  1 &  0 &  1 &  0 &  1
&  0 & \kappa & 3 \\
\hline
   & 1 &  1 &  1 &  1 &  1 &  1 &  1 &  1 &  1 &  1 &  1 &  1 &  1 &  1
&  1 & 1 & 1 \\
D^{\g^7}(Q_8) & 1 &  1 &  0 &  1 &  0 &  1 &  0 &  1 &  0 &  1 &  0
&
1 &  0 &  1 &  0 & 1 & 3 \\
   & 1 &  -1 &  0 &  1 &  0 &  -1 &  0 &  1 &  0 &  -1 &  0 &  1 &  0 &
-1 &  0 & i & 4 \\
   & 2 &  -1 &  0 &  2 &  0 &  -1 &  0 &  2 &  0 &  -1 &  0 &  2 &  0 &
-1 &  0 & 1 & 1 \\
   & 2 &  1 &  0 &  2 &  0 &  1 &  0 &  2 &  0 &  1 &  0 &  2 &  0 &  1
&  0 & -i & 1 \\
   & 2 &  -1 &  0 &  1 &  0 &  -1 &  0 &  0 &  0 &  -1 &  0 &  1 &  0 &
-1 &  0 & \ol\kappa & 3 \\
   & 2 &  -1 &  0 &  1 &  0 &  -1 &  0 &  0 &  0 &  -1 &  0 &  1 &  0 &
-1 &  0 & \kappa^7 & 3 \\
   & 2 &  1 &  0 &  1 &  0 &  1 &  0 &  0 &  0 &  1 &  0 &  1 &  0 &  1
&  0 & \ol\kappa^5 & 3 \\
   & 2 &  1 &  0 &  1 &  0 &  1 &  0 &  0 &  0 &  1 &  0 &  1 &  0 &  1
&  0 & \kappa^3 & 3 \\
\hline
\end{array}
$
\end{tabular}
\newpage
$ \setlength{\arraycolsep}{1pt}
\renewcommand{\arraystretch}{.85}
\begin{array}{|c|c|cccccccccccccc|c|c|}
\hline H &  \chi(1) & \nu^{(2)} & \nu^{(3)}& \nu^{(4)} & \nu^{(5)}&
\nu^{(6)} & \nu^{(7)}& \nu^{(8)} & \nu^{(9)} & \nu^{(10)}&
\nu^{(11)} & \nu^{(12)}& \nu^{(13)} & \nu^{(14)}& \nu^{(15)}& v_\chi
& m \\
\hline
   & 1 &  1 &  1 &  1 &  1 &  1 &  1 &  1 &  1 &  1 &  1 &  1 &  1 &  1
&  1 & 1 & 1 \\
D^{\a_2\a_3}(D_8) & 1 &  1 &  0 &  1 &  0 &  1 &  0 &  1 &  0 &  1 &
0 &  1 &  0 &  1 &  0 & 1 & 3 \\
D^{\a_1\a_3}(D_8) & 1 &  0 &  0 &  1 &  0 &  0 &  0 &  1 &  0 &  0 &
0 &  1 &  0 &  0 &  0 & -i & 4 \\
   & 2 &  -1 &  0 &  1 &  0 &  -1 &  0 &  0 &  0 &  -1 &  0 &  1 &  0 &
-1 &  0 & \ol\kappa^7 & 1 \\
   & 2 &  -1 &  0 &  1 &  0 &  -1 &  0 &  0 &  0 &  -1 &  0 &  1 &  0 &
-1 &  0 & \ol\kappa^3 & 1 \\
   & 2 &  -1 &  0 &  1 &  0 &  -1 &  0 &  0 &  0 &  -1 &  0 &  1 &  0 &
-1 &  0 & \kappa^5 & 1 \\
   & 2 &  -1 &  0 &  1 &  0 &  -1 &  0 &  0 &  0 &  -1 &  0 &  1 &  0 &
-1 &  0 & \kappa & 1 \\
   & 2 &  1 &  0 &  2 &  0 &  1 &  0 &  2 &  0 &  1 &  0 &  2 &  0 &  1
&  0 & 1 & 1 \\
   & 2 &  1 &  0 &  2 &  0 &  1 &  0 &  2 &  0 &  1 &  0 &  2 &  0 &  1
&  0 & i & 1 \\
   & 2 &  0 &  0 &  1 &  0 &  0 &  0 &  0 &  0 &  0 &  0 &  1 &  0 &  0
&  0 & -1 & 2 \\
   & 2 &  0 &  0 &  1 &  0 &  0 &  0 &  0 &  0 &  0 &  0 &  1 &  0 &  0
&  0 & 1 & 2 \\
   & 2 &  0 &  0 &  1 &  0 &  0 &  0 &  0 &  0 &  0 &  0 &  1 &  0 &  0
&  0 & -i & 2 \\
   & 2 &  0 &  0 &  1 &  0 &  0 &  0 &  0 &  0 &  0 &  0 &  1 &  0 &  0
&  0 & i & 2 \\
\hline
   & 1 &  1 &  1 &  1 &  1 &  1 &  1 &  1 &  1 &  1 &  1 &  1 &  1 &  1
&  1 & 1 & 1 \\
D^{\a_1\a_2\a_3}(D_8) & 1 &  1 &  0 &  1 &  0 &  1 &  0 &  1 &  0 &
1 &  0 &  1 &  0 &  1 &  0 & 1 & 3 \\
   & 1 &  0 &  0 &  1 &  0 &  0 &  0 &  1 &  0 &  0 &  0 &  1 &  0 &  0
&  0 & -i & 4 \\
   & 2 &  1 &  0 &  1 &  0 &  1 &  0 &  0 &  0 &  1 &  0 &  1 &  0 &  1
&  0 & \ol\kappa^7 & 1 \\
   & 2 &  1 &  0 &  1 &  0 &  1 &  0 &  0 &  0 &  1 &  0 &  1 &  0 &  1
&  0 & \ol\kappa^3 & 1 \\
   & 2 &  1 &  0 &  1 &  0 &  1 &  0 &  0 &  0 &  1 &  0 &  1 &  0 &  1
&  0 & \kappa^5 & 1 \\
   & 2 &  1 &  0 &  1 &  0 &  1 &  0 &  0 &  0 &  1 &  0 &  1 &  0 &  1
&  0 & \kappa & 1 \\
   & 2 &  1 &  0 &  2 &  0 &  1 &  0 &  2 &  0 &  1 &  0 &  2 &  0 &  1
&  0 & 1 & 1 \\
   & 2 &  1 &  0 &  2 &  0 &  1 &  0 &  2 &  0 &  1 &  0 &  2 &  0 &  1
&  0 & i & 1 \\
   & 2 &  0 &  0 &  1 &  0 &  0 &  0 &  0 &  0 &  0 &  0 &  1 &  0 &  0
&  0 & -1 & 4 \\
   & 2 &  0 &  0 &  1 &  0 &  0 &  0 &  0 &  0 &  0 &  0 &  1 &  0 &  0
&  0 & 1 & 4 \\
\hline
   & 1 &  1 &  1 &  1 &  1 &  1 &  1 &  1 &  1 &  1 &  1 &  1 &  1 &  1
&  1 & 1 & 1 \\
D^{\a_3}(D_8) & 1 &  1 &  0 &  1 &  0 &  1 &  0 &  1 &  0 &  1 &  0
&
1 &  0 &  1 &  0 & 1 & 3 \\
   & 1 &  0 &  0 &  1 &  0 &  0 &  0 &  1 &  0 &  0 &  0 &  1 &  0 &  0
&  0 & -i & 4 \\
   & 2 &  1 &  0 &  1 &  0 &  1 &  0 &  0 &  0 &  1 &  0 &  1 &  0 &  1
&  0 & \ol\kappa^7 & 1 \\
   & 2 &  1 &  0 &  1 &  0 &  1 &  0 &  0 &  0 &  1 &  0 &  1 &  0 &  1
&  0 & \ol\kappa^3 & 1 \\
   & 2 &  1 &  0 &  1 &  0 &  1 &  0 &  0 &  0 &  1 &  0 &  1 &  0 &  1
&  0 & \kappa^5 & 1 \\
   & 2 &  1 &  0 &  1 &  0 &  1 &  0 &  0 &  0 &  1 &  0 &  1 &  0 &  1
&  0 & \kappa & 1 \\
   & 2 &  1 &  0 &  2 &  0 &  1 &  0 &  2 &  0 &  1 &  0 &  2 &  0 &  1
&  0 & 1 & 1 \\
   & 2 &  1 &  0 &  2 &  0 &  1 &  0 &  2 &  0 &  1 &  0 &  2 &  0 &  1
&  0 & i & 1 \\
   & 2 &  0 &  0 &  1 &  0 &  0 &  0 &  0 &  0 &  0 &  0 &  1 &  0 &  0
&  0 & -i & 4 \\
   & 2 &  0 &  0 &  1 &  0 &  0 &  0 &  0 &  0 &  0 &  0 &  1 &  0 &  0
&  0 & i & 4 \\
\hline
   & 1 &  1 &  1 &  1 &  1 &  1 &  1 &  1 &  1 &  1 &  1 &  1 &  1 &  1
&  1 & 1 & 1 \\
D^{\a_2\a_3^3}(D_8) & 1 &  1 &  0 &  1 &  0 &  1 &  0 &  1 &  0 &  1
&  0 &  1 &  0 &  1 &  0 & 1 & 3 \\
D^{\a_1\a_3^3}(D_8) & 1 &  0 &  0 &  1 &  0 &  0 &  0 &  1 &  0 &  0
&  0 &  1 &  0 &  0 &  0 & i & 4 \\
   & 2 &  -1 &  0 &  1 &  0 &  -1 &  0 &  0 &  0 &  -1 &  0 &  1 &  0 &
-1 &  0 & \ol\kappa^5 & 1 \\
   & 2 &  -1 &  0 &  1 &  0 &  -1 &  0 &  0 &  0 &  -1 &  0 &  1 &  0 &
-1 &  0 & \ol\kappa & 1 \\
   & 2 &  -1 &  0 &  1 &  0 &  -1 &  0 &  0 &  0 &  -1 &  0 &  1 &  0 &
-1 &  0 & \kappa^7 & 1 \\
   & 2 &  -1 &  0 &  1 &  0 &  -1 &  0 &  0 &  0 &  -1 &  0 &  1 &  0 &
-1 &  0 & \kappa^3 & 1 \\
   & 2 &  1 &  0 &  2 &  0 &  1 &  0 &  2 &  0 &  1 &  0 &  2 &  0 &  1
&  0 & 1 & 1 \\
   & 2 &  1 &  0 &  2 &  0 &  1 &  0 &  2 &  0 &  1 &  0 &  2 &  0 &  1
&  0 & -i & 1 \\
   & 2 &  0 &  0 &  1 &  0 &  0 &  0 &  0 &  0 &  0 &  0 &  1 &  0 &  0
&  0 & -1 & 2 \\
   & 2 &  0 &  0 &  1 &  0 &  0 &  0 &  0 &  0 &  0 &  0 &  1 &  0 &  0
&  0 & 1 & 2 \\
   & 2 &  0 &  0 &  1 &  0 &  0 &  0 &  0 &  0 &  0 &  0 &  1 &  0 &  0
&  0 & -i & 2 \\
   & 2 &  0 &  0 &  1 &  0 &  0 &  0 &  0 &  0 &  0 &  0 &  1 &  0 &  0
&  0 & i & 2 \\
\hline
\end{array}
$
\newpage
$ \setlength{\arraycolsep}{1pt}
\renewcommand{\arraystretch}{1}
\begin{array}{|c|c|cccccccccccccc|c|c|}
\hline H &  \chi(1) & \nu^{(2)} & \nu^{(3)}& \nu^{(4)} & \nu^{(5)}&
\nu^{(6)} & \nu^{(7)}& \nu^{(8)} & \nu^{(9)} & \nu^{(10)}&
\nu^{(11)} & \nu^{(12)}& \nu^{(13)} & \nu^{(14)}& \nu^{(15)}& v_\chi
& m \\
\hline
   & 1 &  1 &  1 &  1 &  1 &  1 &  1 &  1 &  1 &  1 &  1 &  1 &  1 &  1
&  1 & 1 & 1 \\
D^{\a_1\a_2\a_3^3}(D_8) & 1 &  1 &  0 &  1 &  0 &  1 &  0 &  1 &  0
&
1 &  0 &  1 &  0 &  1 &  0 & 1 & 3 \\
   & 1 &  0 &  0 &  1 &  0 &  0 &  0 &  1 &  0 &  0 &  0 &  1 &  0 &  0
&  0 & i & 4 \\
   & 2 &  1 &  0 &  1 &  0 &  1 &  0 &  0 &  0 &  1 &  0 &  1 &  0 &  1
&  0 & \ol\kappa^5 & 1 \\
   & 2 &  1 &  0 &  1 &  0 &  1 &  0 &  0 &  0 &  1 &  0 &  1 &  0 &  1
&  0 & \ol\kappa & 1 \\
   & 2 &  1 &  0 &  1 &  0 &  1 &  0 &  0 &  0 &  1 &  0 &  1 &  0 &  1
&  0 & \kappa^7 & 1 \\
   & 2 &  1 &  0 &  1 &  0 &  1 &  0 &  0 &  0 &  1 &  0 &  1 &  0 &  1
&  0 & \kappa^3 & 1 \\
   & 2 &  1 &  0 &  2 &  0 &  1 &  0 &  2 &  0 &  1 &  0 &  2 &  0 &  1
&  0 & 1 & 1 \\
   & 2 &  1 &  0 &  2 &  0 &  1 &  0 &  2 &  0 &  1 &  0 &  2 &  0 &  1
&  0 & -i & 1 \\
   & 2 &  0 &  0 &  1 &  0 &  0 &  0 &  0 &  0 &  0 &  0 &  1 &  0 &  0
&  0 & -1 & 4 \\
   & 2 &  0 &  0 &  1 &  0 &  0 &  0 &  0 &  0 &  0 &  0 &  1 &  0 &  0
&  0 & 1 & 4 \\
\hline
   & 1 &  1 &  1 &  1 &  1 &  1 &  1 &  1 &  1 &  1 &  1 &  1 &  1 &  1
&  1 & 1 & 1 \\
D^{\a_3^3}(D_8) & 1 &  1 &  0 &  1 &  0 &  1 &  0 &  1 &  0 &  1 & 0
&  1 &  0 &  1 &  0 & 1 & 3 \\
   & 1 &  0 &  0 &  1 &  0 &  0 &  0 &  1 &  0 &  0 &  0 &  1 &  0 &  0
&  0 & i & 4 \\
   & 2 &  1 &  0 &  1 &  0 &  1 &  0 &  0 &  0 &  1 &  0 &  1 &  0 &  1
&  0 & \ol\kappa^5 & 1 \\
   & 2 &  1 &  0 &  1 &  0 &  1 &  0 &  0 &  0 &  1 &  0 &  1 &  0 &  1
&  0 & \ol\kappa & 1 \\
   & 2 &  1 &  0 &  1 &  0 &  1 &  0 &  0 &  0 &  1 &  0 &  1 &  0 &  1
&  0 & \kappa^7 & 1 \\
   & 2 &  1 &  0 &  1 &  0 &  1 &  0 &  0 &  0 &  1 &  0 &  1 &  0 &  1
&  0 & \kappa^3 & 1 \\
   & 2 &  1 &  0 &  2 &  0 &  1 &  0 &  2 &  0 &  1 &  0 &  2 &  0 &  1
&  0 & 1 & 1 \\
   & 2 &  1 &  0 &  2 &  0 &  1 &  0 &  2 &  0 &  1 &  0 &  2 &  0 &  1
&  0 & -i & 1 \\
   & 2 &  0 &  0 &  1 &  0 &  0 &  0 &  0 &  0 &  0 &  0 &  1 &  0 &  0
&  0 & -i & 4 \\
   & 2 &  0 &  0 &  1 &  0 &  0 &  0 &  0 &  0 &  0 &  0 &  1 &  0 &  0
&  0 & i & 4 \\
\hline
\end{array}
$
\newpage
\begin{tabular}{lcl}
\bf{Frobenius-Schur Exponent 8} & &\bf{Frobenius-Schur Exponent 4}\\
$ \setlength{\arraycolsep}{1pt}
\renewcommand{\arraystretch}{.8}
\begin{array}[t]{|c|c|cccccc|c|c|}
\hline H &  \chi(1) & \nu^{(2)} & \nu^{(3)}& \nu^{(4)} & \nu^{(5)}&
\nu^{(6)} & \nu^{(7)}& v_\chi & m \\
\hline
   & 1 &  1 &  1 &  1 &  1 &  1 &  1 & 1 & 1 \\
D^{\a_2\a_3^2}(D_8) & 1 &  1 &  0 &  1 &  0 &  1 &  0 & 1 & 3 \\
D^{\a_1\a_3^2}(D_8) & 1 &  1 &  0 &  1 &  0 &  1 &  0 & -1 & 4 \\
   & 2 &  -1 &  0 &  0 &  0 &  -1 &  0 & \ol\kappa^6 & 1 \\
   & 2 &  -1 &  0 &  0 &  0 &  -1 &  0 & \ol\kappa^2 & 1 \\
   & 2 &  -1 &  0 &  0 &  0 &  -1 &  0 & \kappa^6 & 1 \\
   & 2 &  -1 &  0 &  0 &  0 &  -1 &  0 & \kappa^2 & 1 \\
   & 2 &  -1 &  0 &  0 &  0 &  -1 &  0 & -i & 2 \\
   & 2 &  -1 &  0 &  0 &  0 &  -1 &  0 & i & 2 \\
   & 2 &  1 &  0 &  0 &  0 &  1 &  0 & -1 & 2 \\
   & 2 &  1 &  0 &  0 &  0 &  1 &  0 & 1 & 2 \\
   & 2 &  1 &  0 &  2 &  0 &  1 &  0 & 1 & 2 \\
\hline
   & 1 &  1 &  1 &  1 &  1 &  1 &  1 & 1 & 1 \\
D^{\a_1\a_2\a_3^2}(D_8) & 1 &  1 &  0 &  1 &  0 &  1 &  0 & 1 & 3 \\
   & 1 &  1 &  0 &  1 &  0 &  1 &  0 & -1 & 4 \\
   & 2 &  1 &  0 &  0 &  0 &  1 &  0 & \ol\kappa^6 & 1 \\
   & 2 &  1 &  0 &  0 &  0 &  1 &  0 & \ol\kappa^2 & 1 \\
   & 2 &  1 &  0 &  0 &  0 &  1 &  0 & \kappa^6 & 1 \\
   & 2 &  1 &  0 &  0 &  0 &  1 &  0 & \kappa^2 & 1 \\
   & 2 &  1 &  0 &  2 &  0 &  1 &  0 & 1 & 2 \\
   & 2 &  -1 &  0 &  0 &  0 &  -1 &  0 & -i & 4 \\
   & 2 &  -1 &  0 &  0 &  0 &  -1 &  0 & i & 4 \\
\hline
   & 1 &  1 &  1 &  1 &  1 &  1 &  1 & 1 & 1 \\
D^{\a_3^2}(D_8) & 1 &  1 &  0 &  1 &  0 &  1 &  0 & 1 & 3 \\
   & 1 &  1 &  0 &  1 &  0 &  1 &  0 & -1 & 4 \\
   & 2 &  1 &  0 &  0 &  0 &  1 &  0 & \ol\kappa^6 & 1 \\
   & 2 &  1 &  0 &  0 &  0 &  1 &  0 & \ol\kappa^2 & 1 \\
   & 2 &  1 &  0 &  0 &  0 &  1 &  0 & \kappa^6 & 1 \\
   & 2 &  1 &  0 &  0 &  0 &  1 &  0 & \kappa^2 & 1 \\
   & 2 &  1 &  0 &  2 &  0 &  1 &  0 & 1 & 2 \\
   & 2 &  1 &  0 &  0 &  0 &  1 &  0 & -1 & 4 \\
   & 2 &  1 &  0 &  0 &  0 &  1 &  0 & 1 & 4 \\
\hline
   & 1 &  1 &  1 &  1 &  1 &  1 &  1 & 1 & 1 \\
D^{\g^6}(Q_8) & 1 &  1 &  0 &  1 &  0 &  1 &  0 & 1 & 3 \\
   & 1 &  1 &  0 &  1 &  0 &  1 &  0 & -1 & 4 \\
   & 2 &  -1 &  0 &  2 &  0 &  -1 &  0 & 1 & 2 \\
   & 2 &  -1 &  0 &  0 &  0 &  -1 &  0 & \ol\kappa^6 & 3 \\
   & 2 &  -1 &  0 &  0 &  0 &  -1 &  0 & \kappa^2 & 3 \\
   & 2 &  1 &  0 &  0 &  0 &  1 &  0 & \ol\kappa^2 & 3 \\
   & 2 &  1 &  0 &  0 &  0 &  1 &  0 & \kappa^6 & 3 \\
\hline
   & 1 &  1 &  1 &  1 &  1 &  1 &  1 & 1 & 1 \\
D^{\g^2}(Q_8) & 1 &  1 &  0 &  1 &  0 &  1 &  0 & 1 & 3 \\
   & 1 &  1 &  0 &  1 &  0 &  1 &  0 & -1 & 4 \\
   & 2 &  -1 &  0 &  2 &  0 &  -1 &  0 & 1 & 2 \\
   & 2 &  -1 &  0 &  0 &  0 &  -1 &  0 & \ol\kappa^2 & 3 \\
   & 2 &  -1 &  0 &  0 &  0 &  -1 &  0 & \kappa^6 & 3 \\
   & 2 &  1 &  0 &  0 &  0 &  1 &  0 & \ol\kappa^6 & 3 \\
   & 2 &  1 &  0 &  0 &  0 &  1 &  0 & \kappa^2 & 3 \\
\hline
\end{array}
$ && $ \setlength{\arraycolsep}{1pt}
\renewcommand{\arraystretch}{.8}
\begin{array}[t]{|c|c|cc|c|c|}
\hline H &  \chi(1) & \nu^{(2)} & \nu^{(3)} & v_\chi & m
\\
\hline
   & 1 &  1 &  1 & 1 & 1 \\
D(Q_8) & 1 &  1 &  0 & 1 & 7 \\
D^{\a_1}(D_8)& 2 &  -1 &  0 & -1 & 1 \\
D^{\a_2}(D_8) & 2 &  -1 &  0 & 1 & 1 \\
D^{\w_{3i}}(E_8) & 2 &  -1 &  0 & -i & 3 \\
   & 2 &  -1 &  0 & i & 3 \\
   & 2 &  1 &  0 & -1 & 3 \\
   & 2 &  1 &  0 & 1 & 3 \\
\hline
   & 1 &  1 &  1 & 1 & 1 \\
D^{\a_1\a_2}(D_8) & 1 &  1 &  0 & 1 & 7 \\
D^{\w_5}(E_8) & 2 &  1 &  0 & -i & 1 \\
   & 2 &  1 &  0 & i & 1 \\
   & 2 &  1 &  0 & -1 & 2 \\
   & 2 &  1 &  0 & 1 & 2 \\
   & 2 &  -1 &  0 & -i & 4 \\
   & 2 &  -1 &  0 & i & 4 \\
\hline
   & 1 &  1 &  1 & 1 & 1 \\
D(D_8) & 1 &  1 &  0 & 1 & 7 \\
D^{\w_1}(E_8) & 2 &  1 &  0 & -i & 1 \\
   & 2 &  1 &  0 & i & 1 \\
   & 2 &  1 &  0 & -1 & 6 \\
   & 2 &  1 &  0 & 1 & 6 \\
\hline
   & 1 &  1 &  1 & 1 & 1 \\
D^{\g^4}(Q_8) & 1 &  1 &  0 & 1 & 7 \\
D^{\w_{3d}}(E_8) & 2 &  1 &  0 & -i & 3 \\
   & 2 &  1 &  0 & i & 3 \\
   & 2 &  -1 &  0 & -1 & 4 \\
   & 2 &  -1 &  0 & 1 & 4 \\
\hline
   & 1 &  1 &  1 & 1 & 1 \\
D^{\w_7}(E_8) & 1 &  1 &  0 & 1 & 7 \\
   & 2 &  1 &  0 & -i & 7 \\
   & 2 &  1 &  0 & i & 7 \\
\hline
\end{array}
$
\end{tabular}


\end{document}